\documentclass[a4paper]{article}
\usepackage[all]{xy}
\usepackage[latin1]{inputenc}        
\usepackage[dvips]{graphics,graphicx}
\usepackage{amsfonts,amssymb,amsmath,color,mathrsfs, amstext,bm,mathtools,tabularx}
\usepackage{amsbsy, amsopn, amscd, amsxtra, amsthm,authblk,enumerate,algorithmicx,algorithm}
\usepackage{booktabs,multirow}
\usepackage{algpseudocode}
\usepackage{upref}
\usepackage{geometry}
\usepackage{lineno}
\usepackage{float}
\usepackage{yhmath}
\usepackage[colorlinks,
            citecolor=red,
            linktocpage]{hyperref}
\usepackage{subfigure}

\DeclareMathOperator{\kl}{KL}
\DeclareMathOperator{\erf}{erf}
\DeclareMathOperator{\erfc}{erfc}

\newcolumntype{b}{>{\hsize=1.5\hsize}X}
\newcolumntype{s}{>{\hsize=.45\hsize}X}

\def\b{\boldsymbol}

\def\e{\varepsilon}

\def\hb {{\hbar}}

\def\bP {\mathbf{P}}

\def\bT {\mathbf{T}}

\def\fH {\mathfrak{H}}
\def\fS {\mathfrak{S}}

\def\cC {\mathcal{C}}
\def\cD {\mathcal{D}}

\def\cF {\mathcal{F}}

\def\cH {\mathcal{H}}

\def\cL {\mathcal{L}}

\def\cO {\mathcal{O}}

\def\a {{\alpha}}
\def\b {{\beta}}
\def\g {{\gamma}}

\def\L {{\Lambda}}
\def\si {{\sigma}}

\def\om {{\omega}}

\def\d {{\partial}}

\def\rstr {{\big |}}
\def\indc {{\bf 1}}

\def\wtilde {\widetilde }

\newcommand{\E}{\mathbb{E}}
\newcommand{\R}{\mathbb{R}}
\newcommand{\bbP}{\mathbb{P}}
\newcommand{\bbZ}{\mathbb{Z}}

\newcommand{\Tr}{\operatorname{trace}}

\newcommand{\ba}{\begin{aligned}}
\newcommand{\ea}{\end{aligned}}

\newcommand{\be}{\begin{equation}}
\newcommand{\ee}{\end{equation}}
\newcommand{\lb}{\label}
\newcommand{\dt}{\Delta t}
\newcommand{\dw}{\Delta W}

\newcommand{\dsp}{\displaystyle}

\newcounter{parentalgorithm}

\makeatother

\newtheorem{theorem}{Theorem}[section]
\newtheorem{lemma}{Lemma}[section]

\newtheorem*{proposition*}{Proposition}

\newtheorem*{corollary*}{Corollary}

\newtheorem*{definitions*}{Definitions}
\newtheorem*{conjecture*}{\bf Conjecture}

\theoremstyle{remark}
\newtheorem{remark}{\bf Remark}[section]

\newtheorem{assumption}{Assumption}[section]

\numberwithin{equation}{section}

\begin{document}

\title{Random Batch Methods for classical and quantum interacting particle systems and statistical samplings}
\author[1]{Shi Jin \thanks{shijin-m@sjtu.edu.cn}}
\author[1]{Lei Li\thanks{leili2010@sjtu.edu.cn}}
\affil[1]{School of Mathematical Sciences, Institute of Natural Sciences, MOE-LSC, Shanghai Jiao Tong University, Shanghai, 200240, P. R. China.}
\date{}
\maketitle

\begin{abstract}
We review the  Random Batch Methods (RBM) for interacting particle systems  consisting of $N$-particles, with $N$ being large. The computational cost of such systems is of $\cO(N^2)$, which is prohibitively expensive.  The RBM  methods use small but random batches so the computational cost is reduced, per time step,  to $\cO(N)$. In this article
we discuss these methods for both classical and quantum systems, the corresponding theory, and applications from molecular dynamics, statistical samplings,   to agent-based models for collective behavior, and quantum Monte-Carlo methods.
\end{abstract}

\tableofcontents

\section{Introduction}

Interacting particle systems arise in a variety of  important phenomena in physical, social, and biological sciences. They usually take the form of  Newton's second law that governs the interactions of $N$-particles under interacting forces that vary depending on
different applications.
Such systems are important in physics--from electrostatics to astrophysics,  in chemistry and  material sciences--such as molecular dynamics,  in biological and social sciences--such as agent based models in  swarming \cite{vicsek1995novel,carrillo2017particle,carlen2013kinetic,degond2017coagulation}, chemotaxis \cite{horstmann03,bertozzi12}, flocking  \cite{cucker2007emergent,hasimple2009,albi2013}, synchronization \cite{choi2011complete,ha2014complete}
 and consensus \cite{motsch2014}).

 These interacting particle systems can be described in general by the first order systems
\begin{gather}\label{eq:Nparticlesys}
d\bm{r}_i=b(\bm{r}_i)\,dt+\alpha_N\sum_{j: j\neq i} K(\bm{r}_i-\bm{r}_j)\,dt+\sigma\, d\bm{W}_i,~~i=1,2,\cdots, N,
\end{gather}
or the second order systems
\begin{gather}\label{eq:Nbody2nd}
\begin{split}
& d\bm{r}_i=\bm{v}_i\,dt,\\
& d\bm{v}_i=\Big[ b(\bm{r}_i)+\alpha_N \sum_{j:j\neq i}K(\bm{r}_i-\bm{r}_j)-\gamma \bm{v}_i \Big]\,dt+\sigma\, d\bm{W}_i.
\end{split}
\end{gather}
We use $\bm{r}_i\in \R^d$ to denote the labels for the particles. 
We will loosely call $\bm{r}_i$ the ``locations'' or ``positions'',
and $\bm{v}_i$ the velocities of the particles, though the specific meaning can be different in different applications. The function $K(\cdot)$ and $b(\cdot): \R^d\to \R^d$ are the interaction kernel and  some given external field respectively. The stochastic processes $\{\bm{W}^i\}_{i=1}^N$
are i.i.d. Wiener processes, or the standard Brownian motions.
 If $\gamma=\sigma=0$ and $b=-\nabla V$ for some potential $V$, one has a Hamiltonian system in classical mechanics.
For the molecules in the heat bath \cite{kawasaki1973simple,callen1951irreversibility}, $\bm{r}_i$ and $\bm{v}_i$ are the physical positions and velocities, described by the underdamped Langevin equations, where $\sigma$ and $\gamma$ satisfy the so-called ``fluctuation-dissipation relation"
\begin{gather}
\sigma=\sqrt{2\gamma/\beta},
\end{gather}
where $\beta$ is the inverse of the temperature (we assume all the quantities are scaled and hence dimensionless so that the Boltzmann constant is absent). The first order system \eqref{eq:Nparticlesys} can be viewed as the overdamped limit (when $\gamma \to \infty$ and the time rescaled as $\gamma t$) of the second order systems \eqref{eq:Nbody2nd}.

In the case $\alpha_N=\frac{1}{N-1}$, as $N \to \infty$, the dynamics of the so-called mean field limit of \eqref{eq:Nparticlesys} is given by the nonlinear Fokker-Planck equation \cite{golse2003mean, mckean1967}
\begin{gather}\label{eq:nonlinearFP}
\partial_t\mu=-\nabla\cdot((b(x)+K*\mu)\mu)+\frac{1}{2}\sigma^2\Delta\mu,
\end{gather}
where $\mu(dx)\in \bP(\R^d)$. The regime $\alpha_N=c/N+o(1)$ is thus naturally called the mean-field regime.
Correspondingly, the mean-field limit of the second order system \eqref{eq:Nbody2nd}  is
\begin{gather}\label{eq:nonlinearLiouville}
\partial_t f=-\nabla_x\cdot(v f)-\nabla_v\cdot((b(x)+K*_x f-\gamma v)f)+\frac{1}{2}\sigma^2\Delta_v f,
\end{gather}
where $f(dx, dv)\in \bP(\R^d\times\R^d)$ and $*_x$ means that the convolution is performed only on the $x$ variable.

If one directly discretizes \eqref{eq:Nparticlesys} or \eqref{eq:Nbody2nd}, the computational cost per time step is $\mathcal{O}(N^2)$. This is undesired for large $N$. The Fast Multipole Method (FMM) \cite{rokhlin1985rapid} is able to reduce the complexity to $\mathcal{O}(N)$ for fast enough decaying interactions. However, the implementation of FMM is quite involved. A simple random algorithm, called the Random Batch Method (RBM), has been proposed in \cite{jin2020random}  to reduce the computation cost per time step from $\mathcal{O}(N^2)$ to $\mathcal{O}(N)$, based on the simple ``random mini-batch" idea. Such an idea is famous for its application in the so-called stochastic gradient descent (SGD) \cite{robbins1951stochastic,bottou1998online,bubeck2015convex} for machine learning problems. The idea was also used for Markov Chain Monte Carlo methods like the stochastic gradient Langevin dynamics (SGLD)  \cite{welling2011bayesian},  and the computation of the mean-field flocking model \cite{albi2013,carrillo2017particle}, motivated by Nanbu's algorithm of the Direct Simulation Monte Carlo method \cite{bird1963approach,nanbu1980direct,babovsky1989convergence}. 

The key behind the ``mini-batch'' idea is to find some cheap unbiased random estimator for the original quantity with the variance being controlled. Depending on the specific applications, the design can be different. For instances, the random batch strategy was proposed regarding general interacting particle systems in \cite{jin2020random}, while the importance sampling in the Fourier space was proposed for the Random Batch Ewald method for molecular dynamics in \cite{jinlixuzhao2020rbe}. Compared with FMM, the accuracy of RBM is lower, but RBM is much simpler and is valid for more general potentials (e.g. the SVGD ODE \cite{li2020svgd}). 
The method converges due to the time average in time, and thus the convergence is like that in the Law of Large Numbers, but in time. 
For long time behaviors,  the method works for systems that own ergodicity and mixing properties, like systems in contact with heat bath and converging to equilibria. A key difference from SGD or SGLD is that the RBM algorithms proposed are aiming to approximate and grasp the dynamical properties of the systems as well, not just to find the optimizer or equilibrium distribution.

RBM for interacting particle systems has been used or extended in various directions,
from statistical sampling \cite{li2020svgd,li2020random,J-L-QMC} to molecular dynamics \cite{jinlixuzhao2020rbe,li2020direct}, control of synchronization \cite{biccari2020stochastic,ko2020model}, and collective behavior of agent-based models \cite{HJKK-C-S-RBM, HJKK3, K-H-J-K3}.
RBM has been shown to converge for finite time interval if the interaction kernels are good enough \cite{li2020svgd,jin2020random}, and in particular an error estimate {\it uniformly} in $N$  was first obtained in  \cite{jin2020random}.  A convergence result of RBM for $N$-body Schr\"odinger equation was established   in \cite{golse2019random}. 

The goal of this review is to introduce the basics of the RBM, the fundamental theory for the convergence and error estimates,
and various applications.

\section{The RBM methods}\label{sec:alg}
In this section, we describe the RBM for general interacting particle systems introduced first
 in \cite{jin2020random}.  We use bold fonts (e.g. $\bm{r}_i$, $\bm{x}_i$, $\bm{v}_i$, $\bm{u}_i$) and capital letters ($X_i, Y_i$)  to denote the quantities that are functions of time $t$ associated with the particles,  use usual letters like $x_i, v_i$ to represent a point in the state space (often $\R^d$), and use letters like $\underline{x}, \underline{v}$ to represent quantities in the configurational space $\R^{Nd}$.

\subsection{The RBM algorithms}

Let $T>0$ be the simulation time, and choose a time step $\dt>0$. Pick a batch size $2\le p\ll N$ that divides $N$ (RBM can also be applied if $p$ does not divide $N$; we assume this only for convenience). Consider the discrete time grids $t_k:=k\dt$, $k\in \mathbb{N}$. For each subinterval $[t_{k-1}, t_k)$, the method has two substeps: (1) at $t_{k-1}$, divide the $N$ particles into $n:=N/p$ groups (batches) randomly; (2) let the particles evolve with interaction only inside the batches.

The above procedure, when applied to the second order system \eqref{eq:Nbody2nd}, leads to Algorithm  \ref{alg:rbm2nd1}. The versions for first order systems is similar.
\begin{algorithm}[H]
\caption{(RBM for \eqref{eq:Nbody2nd})}
\label{alg:rbm2nd1}
\begin{algorithmic}[1]
\For {$m \text{ in } 1: [T/\dt]$}   
\State Divide $\{1, 2, \ldots, N=pn\}$ into $n$ batches randomly.
     \For {each batch  $\mathcal{C}_q$} 
     \State Update $\bm{r}_i, \bm{v}_i$ ($i\in \mathcal{C}_q$) by solving  for $t\in [t_{m-1}, t_m)$ the following
     \begin{gather}\label{eq:RBM2nd}
     \begin{split}
            & d\bm{r}_i=\bm{v}_i\,dt,\\
            & d\bm{v}_i=\Big[b(\bm{r}_i)+\frac{\alpha_N(N-1)}{p-1}\sum_{j\in\mathcal{C}_q,j\neq i}K(\bm{r}_i-\bm{r}_j) -\gamma \bm{v}_i\Big]\,dt+\sigma  d\bm{W}_i.
      \end{split}
      \end{gather}
      \EndFor
 \EndFor
\end{algorithmic}
\end{algorithm}

RBM requires the random division, and the elements in different batches are different.
This is in fact the sampling without replacement. If one allows replacement, one has the following version of RBM.

\begin{algorithm}[H]
\caption{(RBM-r)}\label{randomreplacement}
\begin{algorithmic}[1]
\For {$m \text{ in } 1: [T/\dt]$}   
     \For {$k$ from $1$ to $N/p$}
     \State Pick a set $\mathcal{C}_k$ of size $p$ randomly with replacement. 
     \State Update $\bm{r}_i$'s ($i\in \mathcal{C}_k$) by solving the following SDE for time $\dt$.
     \begin{gather}\label{eq:algorithmreplacement}
     \left\{
           \begin{split}
           & d\bm{x}_i=\bm{u}_i\,dt,\\
            & d\bm{u}_i=\Big[b(\bm{x}_i)+\frac{\alpha_N(N-1)}{p-1}\sum_{j\in\mathcal{C}_q,j\neq i}K(\bm{x}_i-\bm{x}_j) -\gamma \bm{u}_i\Big]\,dt+\sigma\, d\bm{W}_i. \\
            & \bm{x}_i(0) =\bm{r}_i, \quad \bm{u}_i(0)=\bm{v}_i,
            \end{split}
            \right.
      \end{gather}
 \quad\quad\quad i.e., solve \eqref{eq:algorithmreplacement} with initial values $\bm{x}_i(0) =\bm{r}_i, \bm{u}_i(0)=\bm{v}_i$, and set $\bm{r}_i\leftarrow \bm{x}_i(\dt)$, $\bm{v}_i\leftarrow \bm{u}_i(\dt)$.
      \EndFor
\EndFor
\end{algorithmic}
\end{algorithm}

We now discuss the computational cost. Note that random division into $n$ batches of equal size can be implemented using random permutation, which can be realized in $\cO(N)$ operations by Durstenfeld's modern revision of Fisher-Yates shuffle algorithm \cite{durstenfeld1964} (in MATLAB, one can use ``randperm(N)''). After the permutation, one takes the first $p$ elements to be in the first batch, the second $p$ elements to be in the second batch, etc. 
The ODE solver per particle per time step (\ref{eq:algorithmreplacement}) requires merely $\cO(p)$ operations, thus for all particles, each time step costs only $\cO(pN)$. Since $p\ll N$ the overall cost per time step is significantly reduced from  $\cO(N^2)$. 

However, one might encounter the issue of having to use a much smaller time step--which could be of $\cO(N)$ times smaller-- in the RBM implementation. For RBM to really gain significant efficiency, one needs $\dt$ to be {\it independent} of $N$. This is justified by an error analysis to be presented in the next  subsection.

\subsection{Convergence analysis}\label{sec:esterr}

In this subsection, we present the convergence results of RBM for the second order systems \eqref{eq:Nbody2nd} in the mean field regime (i.e., $\alpha_N=1/(N-1)$), which was given   in \cite{jin2020rbm2nd}. 
We remark that the proof relies on the underlying contraction property of the second order systems under certain conditions (\cite{mattingly2002ergodicity,eberle2019couplings}). Due to the degeneracy of the noise terms, the contraction should be proved by suitably chosen variables and Lyapunov functions, and we refer the readers to \cite{jin2020rbm2nd} for more details.

 Denote $(\tilde{\bm{r}}_i, \tilde{\bm{v}}_i)$ the solutions to the RBM process \eqref{eq:RBM2nd}, and use the synchronization coupling as in \cite{jin2020random,jin2021convergence}:
\begin{gather}\label{eq:coupling}
\bm{r}_i(0)=\tilde{\bm{r}}_i(0) \sim \mu_0,~~\bm{W}_i=\tilde{\bm{W}}_i.
\end{gather}

Let $\mathcal{C}_q^{(k)}$ ($1\le q\le n$) be the batches at $t_k$, and  define
\begin{gather}
\mathcal{C}^{(k)}:=\{\mathcal{C}_1^{(k)}, \cdots, \mathcal{C}_n^{(k)}\},
\end{gather}
to be the random division of batches at $t_{k}$.
According to the Kolmogorov extension theorem \cite{durrett2010}, there exists a probability space $(\Omega, \mathcal{F}, \bbP)$ such that the random variables $\{\bm{r}_0^{i}, W^i, \mathcal{C}^{(k)}: 1\le i\le N, k\ge 0\}$ are all defined on this probability space and are independent.
 Let $\E$ denote the integration on $\Omega$ with respect to the probability measure $\bbP$, and consider the $L^2(\cdot)$ norm of a random variable
\begin{gather}
\|\zeta\|=\sqrt{\mathbb{E}|\zeta|^2}.
\end{gather}

For finite time interval, the convergence of RBM is as following. 
\begin{theorem}\label{thm:convfinitetime}
Let $b(\cdot)$ be Lipschitz continuous, and  assume that $|b|, |\nabla b|$ have polynomial growth,  and the interaction kernel $K$ is Lipschitz continuous. Then,
\begin{gather}
\sup_{t\in [0, T]}\sqrt{\E|\tilde{\bm{r}}_i(t)-\bm{r}_i(t)|^2+\E|\tilde{\bm{v}}_i(t)-\bm{v}_i(t)|^2}
\le C(T)\sqrt{\frac{\dt}{p-1}+\dt^2},
\end{gather}
where $C(T)$ is independent of $N$.
\end{theorem}

Often the long-time  error estimates are important since one could use RBM as a sampling method for the invariant measure of \eqref{eq:Nbody2nd} (see section 5). For this we need some additional contraction  assumptions:
\begin{assumption}\label{ass:convexity}
Suppose $b=-\nabla V$ for some $V$ that is bounded from below (i.e., $\inf_x V(x)>-\infty$), and there exist $\lambda_M\ge \lambda_m>0$ such that the eigenvalues of $H:=\nabla^2V$ satisfy
\[
\lambda_m \le \lambda_i(x)\le \lambda_M,~\forall~1\le i\le d, x\in \R^d.
\]
The interaction kernel $K$ is bounded and Lipschitz continuous. Moreover, the friction $\gamma$ and the Lipschitz constant $L$ of $K(\cdot)$ satisfy
\begin{gather}
 \gamma>\sqrt{\lambda_M+2L},~~\lambda_m>2L.
\end{gather}
\end{assumption}
Then the following uniform strong convergence estimate holds:
\begin{theorem}\label{thm:longtimeconv}
Under Assumption \ref{ass:convexity} and the coupling \eqref{eq:coupling}, 
the solutions to \eqref{eq:Nbody2nd} and \eqref{eq:RBM2nd} satisfy
\begin{gather}
\sup_{t\ge 0}\sqrt{\E|\tilde{\bm{r}}_i(t)-\bm{r}_i(t)|^2+\E|\tilde{\bm{v}}_i(t)-\bm{v}_i(t)|^2}
\le C\sqrt{\frac{\dt}{p-1}+\dt^2},
\end{gather}
where the constant $C$ does not depend on $p$ and $N$.
\end{theorem}
Clearly, these error estimates imply that the RBM algorithms can also grasp the dynamical properties. The error estimates above are consequence of some intuitive results, which we summarize here (see \cite{jin2020random}).

For given $\underline{x}:=(x_1, \ldots, x_N)\in \mathbb{R}^{Nd}$, introduce the error of the interacting force
for the $i$th particle.
\begin{gather}
\chi_i(\underline{x}):=\frac{1}{p-1}\sum_{j\in\mathcal{C}}K(x_i-x_j)
-\frac{1}{N-1}\sum_{j:j\neq i}K(x_i-x_j).
\end{gather}
Here, $\mathcal{C}$ is the random batch that contains $i$ in a random division of the batches.

\begin{lemma}\label{lmm:consistency}
Consider a configuration $\underline{x}$ that is independent of the random division. Then,
    \begin{gather}
        \mathbb{E}\chi_i(\underline{x})=0.
    \end{gather}
    Moreover, the (scalar) variance is given by
\begin{gather}\label{eq:var1}
        \mathrm{Var}(\chi_i(\underline{x})) = 
        \left(\frac{1}{p-1}-\frac{1}{N-1}\right)\Lambda_i(\underline{x}),
\end{gather}
    where
\begin{gather}\label{eq:Lambda}
        \Lambda_i(\underline{x}):=\frac{1}{N-2}
        \sum_{j: j\neq i}\Big|  K(x_i-x_j)-\frac{1}{N-1}
        \sum_{\ell: \ell\neq i}K(x_i-x_\ell)  \Big|^2.
\end{gather}
\end{lemma}

Lemma \ref{lmm:consistency} in fact lays the foundation of the convergence of RBM-type
algorithms. The first claim implies that the random estimates of the interacting forces
are unbalanced in the sense that the expectation is zero. This gives the consistency--in expected value--of the RBM approximation,
although each random batch approximation $\frac{1}{p-1}\sum_{j\in\mathcal{C}}K(x_i-x_j)$
to the true interacting force $\frac{1}{N-1}\sum_{j:j\neq i}K(x_i-x_j)$
gives an $\cO(1)$ error (which is clear from $\sqrt{\mathrm{Var}(\chi_i(\underline{x}))}=\cO(1)$). Being a Monte-Carlo like methods, the boundedness of the variance ensures the stability of the RBM methods as can be seen in the proof \cite{jin2020rbm2nd,jin2021convergence}.
 The intrinsic mechanism why such type of methods work is the independent resampling in later time steps, and due to some averaging effect in time
these $\cO(1)$ errors become  small. This Law of Large Numbers type feature {\it in time} guarantees the convergence of RBMs (as indicated by the error bound $\sqrt{|\mathrm{Var}(\chi)|\tau}\sim \sqrt{|\mathrm{Var}(\chi)|/N_T}$ in Theorems \ref{thm:convfinitetime} and \ref{thm:longtimeconv}). 

As another remark, the nonzero variance of the RBM approximation gives some effective noise into the system. This could bring in some ``numerical heating'' effects when RBM is applied for some interacting particle systems. When the system has some dissipation, or in contact with a heat bath as in section section \ref{sec:md} , RBM approximation can be valid for long time and can capture the equilibrium. 

In both Theorem \ref{thm:convfinitetime} and Theorem \ref{thm:longtimeconv}, the error bound is independent of $N$
so that the time step can be chosen independent of $N$ for a fixed accuracy tolerance $\e$.  Hence, for each time step, the cost of RBM is $\cO(1/N)$ of that for direct simulation, but it does not need to take $\cO(N)$ times longer to finish the computation. Such convergence results were first established 
for first order systems \eqref{eq:Nparticlesys} \cite{jin2020random} and then extended to disparate mass cases \cite{jin2021convergence}. The weak convergence has also been discussed in \cite{jin2021convergence}.

\subsection{An illustrating example: wealth evolution}

To illustrate the algorithms, we consider the model proposed by Degond et. al. \cite{degond2014} for
the evolution of $N$ market agents with two attributes: the economic configuration $X_i$ and its wealth $Y_i$.
\begin{gather}\label{eq:ecowealth}
\begin{split}
& \dot{X}_i=V(X_i, Y_i),\\
& dY_i=-\frac{1}{N-1}\sum_{k:k\neq i}\xi_{ik}\Psi(|X_i-X_k|)\partial_y\phi(Y_i-Y_k)\,dt+\sqrt{2D} Y_i dW_i.
\end{split}
\end{gather}
The first equation describes the evolution of the economic configuration, which is driven by the local Nash equilibrium and it is related to mean-field games \cite{lasry2007}. This model clearly is a suitable interacting particle system for which the RBM algorithms fit perfectly. Moreover, the RBM version of \eqref{eq:ecowealth} can be viewed as a new model as one agent may only trade with a small number of random agents during a short time in the real world.

For numerical experiments, \cite{jin2020random} considers the homogeneous case when the wealth dynamics is independent of the position in the economic configuration space. Then, the dynamics of the wealth is reduced the following
\begin{gather}\label{eq:wealth}
 dY_i=-\frac{\kappa}{N-1} \sum_{k:k\neq i}\partial_y\phi(Y_i-Y_k)\,dt+\sqrt{2D} Y_i dW_i.
\end{gather}
The corresponding mean field dynamics has an equilibrium distribution given by
\[
\rho_{\infty}(y)\propto \exp\left(-\frac{\alpha(y)}{D}\right),
\]
where $\alpha$ satisfies
\[
\partial_y\alpha(y)=-\frac{1}{y^2}F(y)+\frac{2D}{y}.
\]

\begin{figure}
\begin{center}
	\includegraphics[width=0.8\textwidth]{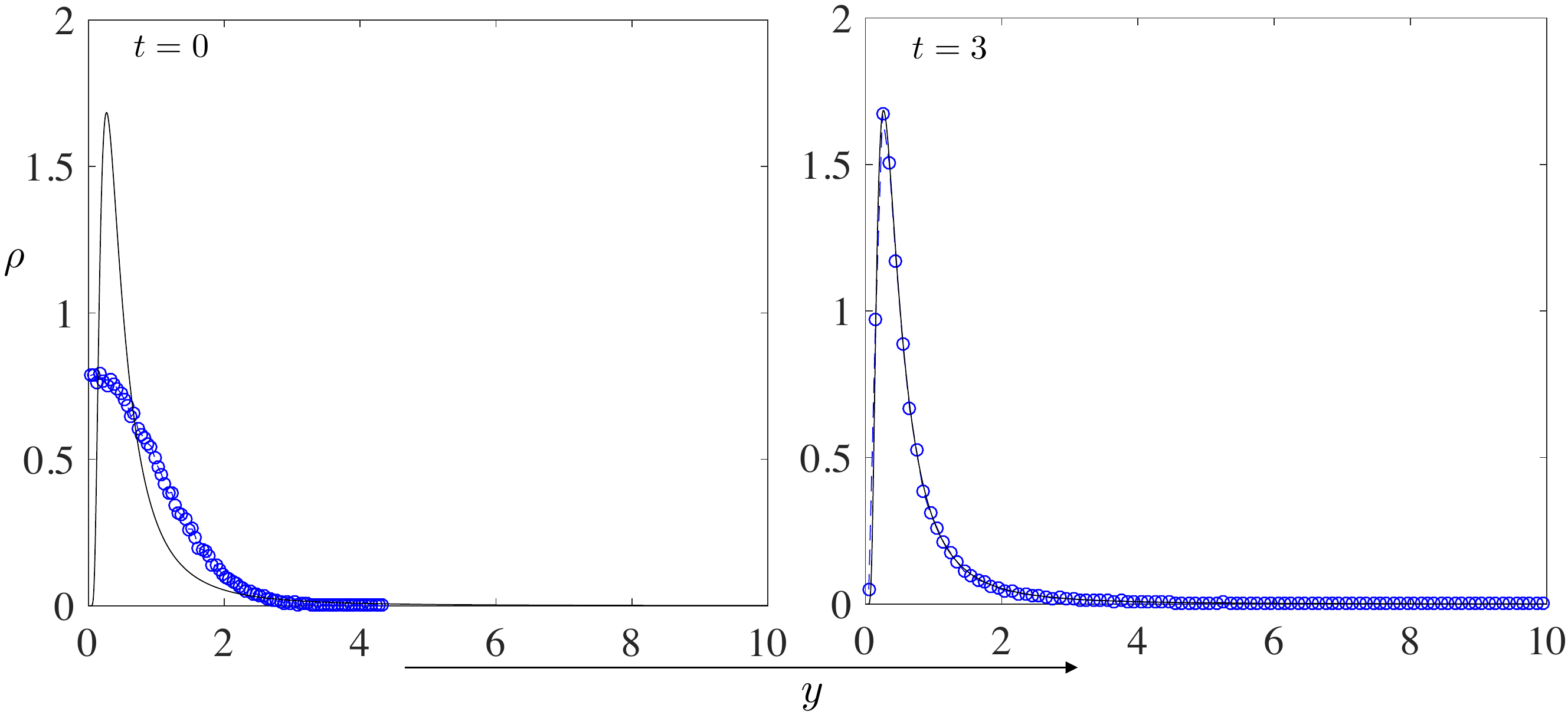}
\end{center}
\caption{Wealth distribution obtained by RBM compared with the reference curve}
\label{fig:wealth}
\end{figure}

In Fig. \ref{fig:wealth}, the empirical distribution of the wealth obtained by RBM for the case $\phi(y)=\frac{1}{2}y^2$ is compared to the reference curve (an inverse Gamma distribution), which is 
\[
\rho_{\infty}(y)=\frac{(\kappa \eta/D)^{\kappa/D+1}}{\Gamma(\kappa/D+1)}y^{-(2+\kappa/D)}\exp\left(-\frac{\kappa \eta}{Dy}\right)1_{y>0},~~\eta=\frac{\sqrt{2}}{\sqrt{\pi}}.
\]
Clearly, the distribution obtained by RBM agrees perfectly with the expected wealth distribution at $t=3$ already.

This example has two distinguished features: long range and multiplicative noises. We point out that although the convergence results presented in subsection \ref{sec:esterr} were for regular interacting potentials $K$ and additive noises, as shown by this and more examples in \cite{jin2020random}, and those in later sections,  the RBM algorithms are applicable to much broader classes of interacting particle systems, including long-range, singular (like the Lenard-Jones and Coulomb) potentials (see section \ref{sec:md} below), and with multiplicative noise.

\section{The mean-field limit}

It is known that the $N$-particle system \eqref{eq:Nparticlesys} has the mean-field limit given by the Fokker-Planck equation
\eqref{eq:nonlinearFP}. Namely, the empirical measure of the particle system \eqref{eq:Nparticlesys} is close, in Wasserstein distance, to $\mu$ in \eqref{eq:nonlinearFP}. Thus, when $N$ is large, one may use the RBM as a numerical (particle method) for 
\eqref{eq:nonlinearFP}. Indeed, since the error bounds obtained 
in the precious section is independent of $N$, one could hope
that when $N \to \infty$, the empirical measure of the RBM 
should be close to $\mu$. To justify this, one first needs to
derive the mean field limit of the RBM, for fixed $\Delta t$, then compare it with
\eqref{eq:nonlinearFP}. In addition, the RBM could be viewed
as a random model for the underlying physics, hence it is
also natural to ask what its mean field limit is.

 The mean-field limit of the RBM for the general first order system \eqref{eq:Nparticlesys}, given below by Algorithm \ref{rbm_first}, was derived and proved in  \cite{jin2021mean}. We summarize the results in this section. 

\begin{algorithm}[H]
\caption{(RBM for first order systems)}\label{rbm_first}
\begin{algorithmic}[1]
\For {$k \text{ in } 1: [T/\dt]$}   
\State Divide $\{1, 2, \ldots, N\}$ into $n=N/p$ batches randomly.
     \For {each batch  $\mathcal{C}_q$} 
     \State Update $\bm{r}_i$'s ($i\in \mathcal{C}_q$) by solving the following SDE with $t\in [t_{k-1}, t_k)$.
     \begin{gather}\label{eq:rbmSDE}
            d \bm{r}_i=b(\bm{r}_i) dt+\frac{1}{p-1}\sum_{j\in \mathcal{C}_q,j\neq i} K(\bm{r}_i-\bm{r}_j)dt+\sigma\, d\bm{W}^i.
      \end{gather}
      \EndFor
\EndFor
\end{algorithmic}
\end{algorithm}

Intuitively, when $N\gg 1$, the probability that two chosen particles are correlated is very small. Hence, in the $N\to\infty$ limit, two chosen particles will be independent with probability $1$. Due to the exchangeability, the marginal distributions of the particles will be identical.  Based on this observation, the following mean field limit for RBM can be obtained for the one-particle distribution: 
\begin{algorithm}[H]
\caption{(Mean Field Dynamics of RBM \eqref{eq:rbmSDE})}\label{meanfield}
\begin{algorithmic}[1]
\State $\tilde{\mu}(\cdot, t_0)=\mu_0$.
\For {$k \ge 0$}  

\State Let $\rho^{(p)}(\cdots, 0)=\tilde{\mu}(\cdot, t_{k})^{\otimes p}$ be a probability measure on $(\R^{d})^{ p}\cong \R^{pd}$.

\State Evolve the measure $\rho^{(p)}$ to find $\rho^{(p)}(\cdots, \dt)$ by the following Fokker-Planck equation:
\begin{gather}\label{eq:firstalgorithm}
            \partial_t\rho^{(p)}=-\sum_{i=1}^p 
            \nabla_{x_i}\cdot\left(\Big[b(x_i)+\frac{1}{p-1}\sum_{j=1,j\neq i}^p K(x_i-x_j)\Big]\rho^{(p)}\right)+\frac{1}{2}\sigma^2\sum_{i=1}^p \Delta_{x_i}\rho^{(p)}.
\end{gather}

\State Set
\begin{gather}
\tilde{\mu}(\cdot, t_{k+1}):=\int_{(\R^{d})^{\otimes(p-1)}}
\rho^{(p)}(\cdot,dy_2,\cdots,dy_p, \dt).
\end{gather}

\EndFor
\end{algorithmic}
\end{algorithm}
The dynamics in Algorithm \ref{meanfield} naturally gives a nonlinear operator $\mathcal{G}_{\infty}: \bP(\R^d)\to \bP(\R^d)$ as
\begin{gather}\label{eq:Ginfty}
\tilde{\mu}(\cdot, t_{k+1})=: \mathcal{G}_{\infty}(\tilde{\mu}(\cdot, t_k)).
\end{gather}
Corresponding to this is the following SDE system for $t\in [t_k, t_{k+1})$
\begin{gather}\label{eq:RBMmeanfieldSDE}
d\bm{x}_i=b(\bm{x}_i)\,dt+\frac{1}{p-1}\sum_{j=1,j\neq i}^{p}K(\bm{x}_i-\bm{x}_j)\,dt
+\sigma\,d\bm{W}_i,~~i=1,\cdots, p,
\end{gather}
with $\{\bm{x}_i(t_k)\}$ drawn i.i.d from $\tilde{\mu}(\cdot, t_k)$.  Then, $\tilde{\mu}(\cdot, t_{k+1})=\mathscr{L}(\bm{x}_1(t_{k+1}^-))$, the law of $\bm{x}_1(t_{k+1}^-)$. Note that all $\bm{x}_i$ have the same distribution for any $t_k\le t< t_{k+1}$.
Without loss of generality, we will impose $\bm{x}_1(t_k^-)=\bm{x}_1(t_k^+)$. For other particles $i\neq 1$,  $\bm{x}_i(t)$ in $[t_{k-1}, t_k)$ and $[t_k, t_{k+1})$ are independent and they are not continuous at $t_k$. In fact, in the $N\to\infty$ limit, $\bm{x}_i, i\neq 1$ at different subintervals correspond to different particles that interact with particle $1$ as in Algorithm \ref{rbm_first}.

Hence, in the mean field limit of RBM, one starts with a chaotic configuration\footnote{By "chaotic configuration", we mean that there exists a one particle distribution $f$ such that for any $j$, the $j$-marginal distribution is given by $\mu^{(j)}=f^{\otimes j}$. Such independence in a configuration is then loosely called "chaos". If the $j$-marginal distribution is more close to $f^{\otimes j}$ for some $f$, we loosely say "there is more chaos".}, the $p$ particles evolve by interacting with  each other. Then, at the starting point of the next time interval, one imposes the chaos condition so that the particles are independent again. 

In \cite{jin2021mean}, this intuition has been justified rigorously for finitely many steps under the following assumptions. 

\begin{assumption}\label{ass:momentmu0}
The moments of the initial data are finite:
\begin{gather}
\int_{\R^d} |x|^q\mu_0(dx)<\infty,~\forall q\in [2,\infty).
\end{gather}
\end{assumption}

\begin{assumption}\label{ass:kernelfunctions}
Assume $b(\cdot): \R^d\to\R^d$ and $K(\cdot): \R^d\to \R^d$ satisfy the following conditions.
\begin{itemize}
\item  $b(\cdot)$ is one-sided Lipschitz:
\begin{gather}
(z_1-z_2)\cdot (b(z_1)-b(z_2))\le \beta |z_1-z_2|^2
\end{gather}
for some constant $\beta$;

\item $K$ is Lipschitz continuous
\[
|K(z_1)-K(z_2)|\le L|z_1-z_2|.
\]
\end{itemize}
\end{assumption}

\begin{figure}
\begin{center}
	\includegraphics[width=0.8\textwidth]{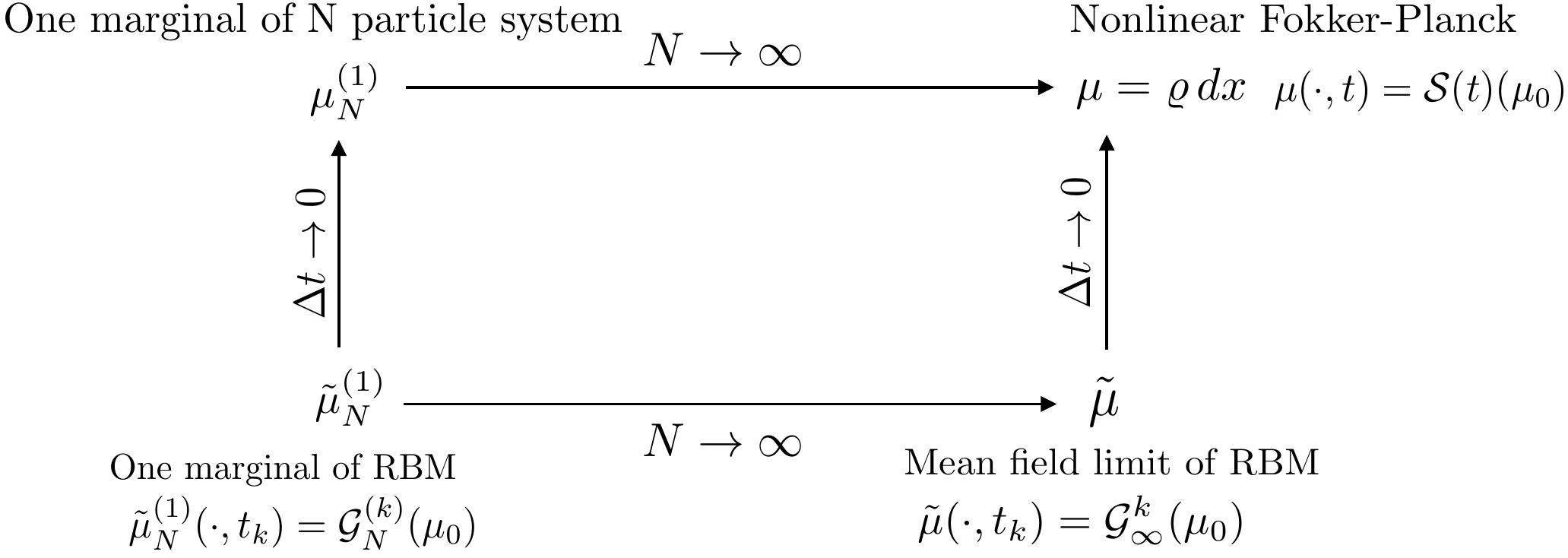}
\end{center}
\caption{Illustration of the various operators and the asymptotic limits.}
\label{fig:operator}
\end{figure}

Corresponding to the operator \eqref{eq:Ginfty}, one may define the operator $\mathcal{G}_N^{k}: \bP(\R^d)\to \bP(\R^d)$
for RBM with $N$ particles as follows.
Let $\bm{r}_i(0)$'s be i.i.d drawn from $\mu_0^{\otimes N}$, and consider \eqref{eq:rbmSDE}. Define
\begin{gather}
\mathcal{G}_N^{k}(\mu_0):=\mathscr{L}(\bm{r}_1(t_k)).
\end{gather}
Recall that $\mathscr{L}(\bm{r}_1)$ denotes the law of $\bm{r}_1$, thus the one marginal distribution. Conditioning on a specific sequence of random batches, the particles are not exchangeable. However,
when one considers the mixture of all possible sequences of random batches, the laws of the particles $\bm{r}_i(t_k)$ ($1\le i\le N$) are identical. In Fig.\ref{fig:operator}, we illustrate these definitions and various limits.
With these setup introduced, we may state the first main result in \cite{jin2021mean} as follows:
\begin{theorem}\label{thm:meanfield}
Under assumptions \ref{ass:momentmu0} and \ref{ass:kernelfunctions},  for any fixed $k$, it holds that
\begin{gather}
\lim_{N\to\infty}W_q(\mathcal{G}_{\infty}^k(\mu_0), \mathcal{G}_N^{k}(\mu_0))=0.
\end{gather}
\end{theorem}
Here, $W_q$ is the Wasserstein-q distance  \cite{santambrogio2015}:
\begin{gather}\label{eq:Wassersteindef}
W_q(\mu, \nu)=\left(\inf_{\gamma\in \Pi(\mu,\nu)}\int_{\R^d\times \R^d} |x-y|^q d\gamma\right)^{1/q},
\end{gather}
where $\Pi(\mu, \nu)$ is the set of ``transport plans", i.e. a joint measure on $\R^d\times \R^d$ such that the marginal measures corresponding to $x$ and $y$ are $\mu$ and $\nu$ respectively.

The next questions is whether 
 the one marginal distribution $\mu_N^{(1)}:=\mathscr{L}(\bm{r}_1)$ of the RBM converges to $\mu$. Denote the solution operator to \eqref{eq:nonlinearFP} by $\mathcal{S}$:
\begin{gather}
\mathcal{S}(\Delta)\mu(t_1):=\mu(t_1+\Delta),~\forall t_1\ge 0, \Delta\ge 0.
\end{gather}
Clearly, $\{\mathcal{S}(t): t\ge 0\}$ is a nonlinear semigroup.

We make  more technical assumptions here.
\begin{assumption}\label{ass:newrho0}
The measure $\mu_0$ has a density $\varrho_0$ that is smooth with finite moments
$\int_{\R^d}|x|^q\varrho_0\,dx<\infty$, $\forall q\ge 1$,
and the entropy is finite
\begin{gather}
H(\mu_0):=\int_{\R^d}\varrho_0\log\varrho_0\,dx<\infty.
\end{gather}
\end{assumption}
If $\varrho_0(x)=0$ at some point $x$, one defines $\varrho_0(x)\log\varrho_0(x)=0$.
We also introduce the following assumption on the growth rate of derivatives of $b$ and $K$, which will be used below.
\begin{assumption}\label{ass:polynomialgrowth}
The function $b$ and its derivatives have polynomial growth.
The derivatives of $K$ with order at least $2$ (i.e., $D^{\alpha}K$ with $|\alpha|\ge 2$) have polynomial growth.
\end{assumption}

Based on these conditions, it can be shown that $\mu$ has a density $\varrho(\cdot, t)$. For convenience, 
we will not distinguish $\mu$ from its density $\varrho$. Sometimes, one may also assume the strong confinement condition:
\begin{assumption}\label{ass:kernelstrong}
The fields $b(\cdot): \R^d\to\R^d$ and $K(\cdot): \R^d\to \R^d$ are smooth. Moreover,  $b(\cdot)$ is strongly confining:
\begin{gather}
(z_1-z_2)\cdot (b(z_1)-b(z_2))\le -r |z_1-z_2|^2
\end{gather}
for some constant $r>0$, and $K$ is Lipschitz continuous
$|K(z_1)-K(z_2)|\le L|z_1-z_2|$.
The parameters $r, L$ satisfy
\begin{gather}
r>2L.
\end{gather}
\end{assumption}

With the assumptions stated, we can state the second main result in \cite{jin2021mean}.
\begin{theorem}\label{thm:w2distance}
Suppose Assumptions \ref{ass:kernelfunctions}, \ref{ass:newrho0} and \ref{ass:polynomialgrowth} hold. Then, 
\begin{gather}
\sup_{n: n\dt\le T}W_{1}(\mathcal{G}_{\infty}^n(\varrho_0), \varrho(n\dt))\le C(T)\dt.
\end{gather}

If Assumption \ref{ass:kernelstrong} is assumed in place of Assumption \ref{ass:kernelfunctions} and also $\sigma>0$, then
\begin{gather}
\sup_{n\ge 0}W_{1}(\mathcal{G}_{\infty}^n(\varrho_0), \varrho(n\dt))\le C\dt.
\end{gather}
\end{theorem}

 These theorems show that the dynamics given by $\mathcal{G}_{\infty}$ can approximate that of the nonlinear Fokker-Planck equation \eqref{eq:nonlinearFP}, with the $W_1$ distance to be of $\cO(\dt)$ . Thus, the two limits $\lim_{N\to \infty}$ and $\lim_{\dt\to 0}$ commute, as shown in Fig. \ref{fig:operator}.

\section{Molecular dynamics}\label{sec:md}

Molecular dynamics (MD) refers to computer simulation of atoms and molecules, and is among the most popular numerical methods to understand the dynamical and equilibrium properties of many-body particle systems in many areas such as chemical physics, soft materials and biophysics \cite{ciccotti1987simulation,frenkel2001understanding,FPP+:RMP:2010}.   In this section, we discuss the relevant issues and the applications of RBM and its modifications.

Consider $N$ ``molecules'' with  masses $m_i$'s (each might be a model for a real molecule or a numerical molecule that is a packet of many real molecules) that interact with each other.  The equations of motion are given by
\begin{gather}\label{eq:md1}
\begin{split}
& d\bm{r}_i=\bm{v}_i\,dt,\\
& m_i d\bm{v}_i=\Big[-\sum_{j:j\neq i}\nabla \phi(\bm{r}_i-\bm{r}_j)\Big]\,dt
+d\bm{\xi}_i.
\end{split}
\end{gather}
Here, $\phi(\cdot)$ is the interaction potential and $d\bm{\xi}_i$ means some other possible terms that change the momentum,. Typical examples of the potential include the Coulomb potentials 
\[
\phi(x) = \frac{q_i q_j}{r},
\]
where $q_i$ is the charge for the $i$th particle and $r=|x|$, and the Lennard-Jones potential 
\[
\phi(x)=4\left(\frac{1}{r^{12}}-\frac{1}{r^6}\right).
\]
Between ions, both types of potential exist and between charge-neutral molecules, the Lennard-Jones potential might be the main force (the Lennard-Jones interaction intrinsically also arises from the interactions between charges, so these two types are in fact both electromagnetic forces) \cite{frenkel2001understanding,FPP+:RMP:2010}.  To model the solids or fluids with large volume, one often uses a box with length $L$, equipped with the periodic conditions for the simulations. 

To model the interaction between the molecules with the heat bath, one may consider some thermostats so that the temperature of the system can be controlled at a given value. The thermostats are especially good for RBM approximations as the effective noise introduced by RBM approximation can be damped by the thermostats, reducing the ``numerical heating'' effects \cite{jinlixuzhao2020rbe}.  Typical thermostats include the Andersen thermostat, the Langevin thermostat and the Nos\'e-Hoover thermostat \cite{frenkel2001understanding}. 
In the Andersen thermostat \cite[section 6.1.1]{frenkel2001understanding}, one does the simulation for 
$d\bm{\xi}_i=0$ 
between two time steps, but a particle can collide with the heat bath at each discrete time. Specifically, assume the collision frequency is $\nu$, so in a duration of time $t\ll 1$ the chance that a collision has happened is given by the exponential distribution
\[
1-\exp(-\nu t)\approx \nu t,~~t\ll 1.
\]
If a collision happens, the new velocity is then sampled from the Maxwellian distribution with temperature $T$ (i.e., the normal distribution $\mathcal{N}(0, T)$). In the underdamped Langevin dynamics, one chooses 
\[
d\bm{\xi}_i=-\gamma \bm{v}_i\,dt+\sqrt{\frac{2\gamma}{\beta}}\,d\bm{W}_i,
\] 
so that the ``fluctuation-dissipation relation'' is satisfied and the system will evolve to the equilibrium with the correct temperature $T=\beta^{-1}$. 
It is well-known that the invariant measure of such systems is given by the  Gibbs distribution \cite{lifshitz2013statistical}
\[
\pi(\underline{x}, \underline{v}) \propto \exp\left(-\beta(\frac{1}{2}\sum_{i=1}^N |v_i|^2+U(\underline{x}))\right),
U(\underline{x})=\frac{1}{2}\sum_{i,j: i\neq j}\phi(x_i-x_j),
\]
where $\underline{x}=(x_1,\cdots, x_N)\in \R^{Nd}$ and $\underline{v}=(v_1,\cdots, v_N)\in \R^{Nd}$.
The Nos\'e-Hoover thermostat uses a Hamiltonian for an extended system of $N$ particles plus an additional coordinate $s$ (\cite{nose1984molecular,hoover1985canonical}):
\[
\mathcal{H}_{\mathrm{NH}}=\sum_{i=1}^N\frac{|\tilde{\bm{p}}_i|^2}{2 m_i s^2}+U(\{\bm{r}_i\})
+\frac{p_s^2}{2Q}+L\frac{\ln s}{\beta}.
\]
Here, $\tilde{\bm{p}}_i$ is the momentum of the $i$th particle.
The microcanonical ensemble corresponding to this Hamiltonian reduces to the canonical ensemble for the real variables $\bm{p}_i=\tilde{\bm{p}}_i/s$. Hence, one may run the following deterministic ODEs, which are the Hamiltonian ODEs with Hamiltonian $\mathcal{H}_{\mathrm{NH}}$ in terms of the so-called real variables,
\[
\begin{split}
&\dot{\bm{r}}_i=\bm{p}_i,\\
& \dot{\bm{p}}_i=-\nabla_{\bm{r}_i}U-\xi \bm{p}_i,\\
& \dot{\xi}=\frac{1}{Q}\left(\sum_{i=1}\frac{|\bm{p}_i|^2}{m_i}-\frac{3N}{\beta}\right).
\end{split}
\]
The time average of the desired quantities will be the correct canonical ensemble average. As one can see, when the temperature of the system is different from $T$,
the extra term $-\xi \bm{p}_i$  will drive the system back to temperature $T$,  thus it may give better behaviors for controlling the temperature.

\subsection{RBM with kernel splitting}\label{subsec:splitting}

In molecular dynamics simulation, the interaction force kernel 
\[
K(x):=-\nabla\phi(x),~x\in\R^d,
\]
is often singular at $x=0$. Hence, the direct application of RBM could lead to poor results. To resolve this issue, one can adopt the splitting strategy in \cite{martin1998novel,hetenyi2002multiple}, and decompose the interacting force $K$ into two parts:
\begin{gather}
K(x)=K_1(x)+K_2(x).
\end{gather}
Here, $K_1$ has short range that vanishes for $|x|\ge r_0$ where $r_0$ is a certain cutoff chosen to be comparable to the mean distance of the particles. $K_2(x)$ is a bounded smooth function.
One then applies RBM to the $K_2$ part only. The resulted method is shown in Algorithm \ref{alg:rbms}.
Now, the cost of summation in $K_1$ is of $\mathcal{O}(1)$  for each given $i$ using data structures like Cell-List \cite[Appendix F]{frenkel2001understanding}. Since $K_2$ is bounded, RBM can be applied well due to the boundedness of variance, without introducing too much error. 
Hence, the cost per time step is again $\mathcal{O}(N)$.  For practical applications, one places the initial positions of the molecules on the grid of a lattice, and the repulsive force $K_1$ will forbid the particles from getting too close so that the system is not too stiff. 

\begin{algorithm}[H]
\caption{RBM with splitting for \eqref{eq:Nbody2nd} }
\label{alg:rbms}
\begin{algorithmic}[1]
\State Split $K=:K_1+K_2$, where $K_1$ has short range, while $K_2$ has long range but is smooth.

\For {$m \text{ in } 1: [T/\dt]$}   
\State Divide $\{1, 2, \ldots, N=pn\}$ into $n$ batches randomly.
     \For {each batch  $\mathcal{C}_q$} 
     \State Update $(\bm{r}_i, \bm{v}_i)$'s ($i\in \mathcal{C}_q$) by solving for $t\in [t_{m-1}, t_m)$
     \begin{gather}\label{eq:rbm2ndorder}
     \begin{split}
             d\bm{r}_i=\,&\bm{v}_i\,dt ,\\
             d\bm{v}_i=\,&\Big[ b(\bm{r}_i)+\alpha_N\sum_{j: j\neq i}K_1(\bm{r}_i-\bm{r}_j)-\gamma \bm{v}_i\Big]\,dt\\
            &+\frac{\alpha_N(N-1)}{p-1}\sum_{j\in\mathcal{C}_q,j\neq i}K_2(\bm{r}_i-\bm{r}_j)\,dt+\sigma\, d\bm{W}_i.
      \end{split}
      \end{gather}
      \EndFor
 \EndFor
\end{algorithmic}
\end{algorithm}

\begin{figure}[!htbp]
\centering
\subfigure[$\nu=\gamma=10$, $\dt_k=0.001/\log(k+1)$.]{
\begin{minipage}[b]{0.47\textwidth}
\includegraphics[width=1\textwidth]{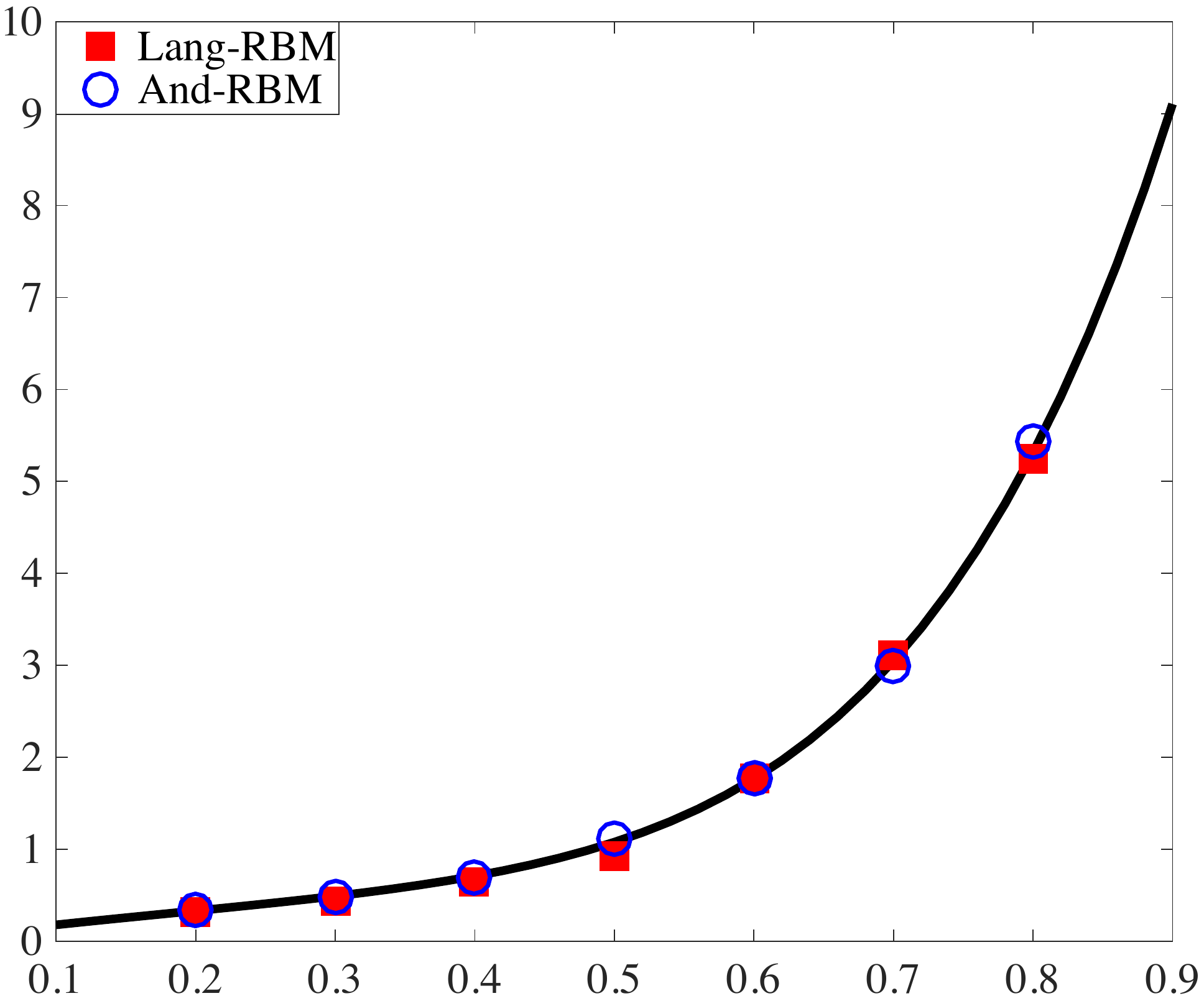}
\end{minipage}}
\subfigure[$\nu=\gamma=50$, $\dt=0.001$.]{
\begin{minipage}[b]{0.47\textwidth}
\includegraphics[width=1\textwidth]{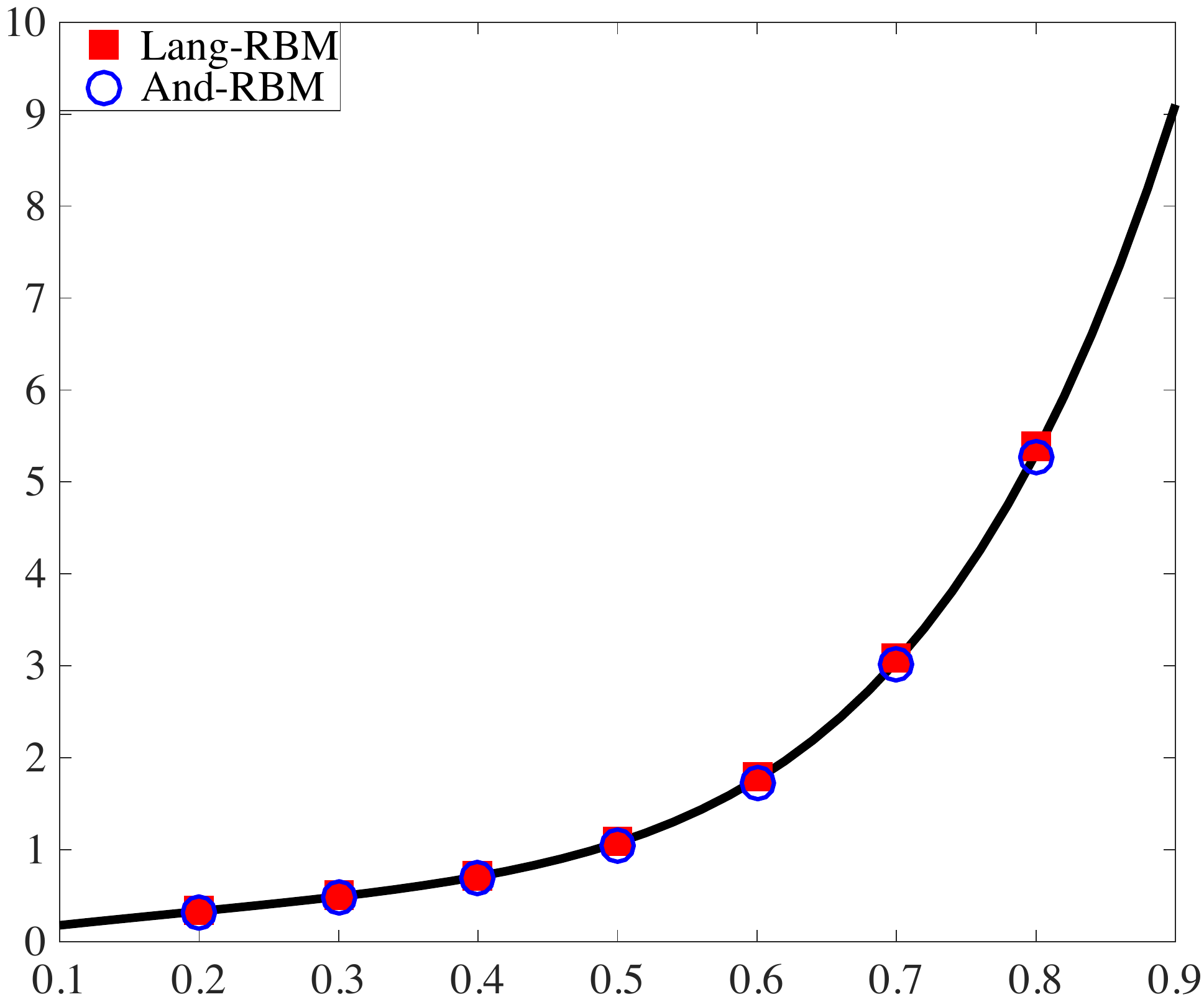}
\end{minipage}
}
\caption{The pressure obtained by Andersen-RBM and Langevin-RBM for Lenard-Jones fluid with $N=500$: the blue circles are those by Andersen-RBM while the red squares are by Langevin-RBM.}
\label{fig:lgvadr}
\end{figure}

Using this splitting strategy, one may apply RBM  to the MD simulations with different thermostats.
In Fig. \ref{fig:lgvadr}, we show the numerical results from \cite{jin2020rbm2nd} for a Lennard-Jones fluid with temperature $\beta^{-1}=2$
and the length of box $L=(N/\rho)^{1/3}$ for a given density $\rho$. The results are obtained using the Andersen thermostat and the Langevin thermostat respectively, with the splitting and RBM strategy, for particle number $N=500$. 
The reference curves (black solid line) are the fitting curves in \cite{johnson1993lennard}. 
In the first figure, the decreasing step sizes $\dt_k=0.001/\log(k+1)$ are taken to reduce the numerical heating effect brought by RBM
when the collision coefficient are not so big ($\nu=\gamma=10$).
The results show that RBM with splitting strategy can work reasonably well for the Lennard-Jones fluid in the considered regime.

\subsection{Random Batch Ewald: an importance sampling in the Fourier space}

In the presence of long-range interactions such as the Coulomb interactions, the molecular dynamics simulations becomes computationally expensive for large $N$. A lot of effort has already  been  devoted to computing such long-range interactions efficiently. Some popular methods include lattice summation methods such as the particle-particle particle mesh Ewald (PPPM) \cite{LDT+:MS:1994,Deserno98JCP}, and multipole type methods such as treecode \cite{BH:N:1986,DK:JCP:2000} and fast multipole methods (FMM) \cite{GR:JCP:1987,YBZ:JCP:2004}. 
These methods can reduce the complexity per time step from $O(N^2)$ to $\cO(N\log N)$ or even $\cO(N)$, and have gained big success in practice. However, some issues still remain to be resolved, e.g., the prefactor in the linear scaling can be large, or the implementation can be nontrivial, or the scalability for parallel computing is not high.

In this section, we give a brief introduction to the recently proposed  Random Batch Ewald (RBE) method  for molecular dynamics simulations of particle systems with long-range Coulomb interactions, which achieves an $\cO(N)$ complexity \cite{jinlixuzhao2020rbe} with a high parallel efficiency.  The RBE method is based on the Ewald splitting for the Coulomb kernel with a random ``mini-batch'' type technique applied in the Fourier series for the long-range part. 
Compared with PPPM where the Fast Fourier Transform is used to speed up the computation in the Fourier space, the RBE method uses random batch type technique to speed up the computation.

Consider $N$ physical or numerical particles
inside the periodic box with side length $L$, assumed to have net charge $q_i$ ($1\le i\le N$) with the electroneutrality condition
\begin{gather}
\sum_{i=1}^N q_i=0.
\end{gather}
The forces are computed using $\bm{F}_i=-\nabla_{\bm{r}_i}U$, where $U$ is the potential energy of the system.
Since the Coulomb potential is of long range, with the periodic boundary condition, one must consider the images so that
\begin{gather}\label{eq:energy}
U=\frac{1}{2}\sum_{\bm{n}}{}'\sum_{i,j=1}^N q_iq_j \frac{1}{|\bm{r}_{ij}+\bm{n}L|},
\end{gather}
where $\bm{n}\in \mathbb{Z}^3$ ranges over the three-dimensional integer vectors and $\sum'$ is defined such that $\bm{n}=0$ is not included when $i=j$.

Due to the long-range nature of the Coulomb potential, the series \eqref{eq:energy} converges conditionally, thus a naive truncation would require  a very large  $r$ to maintain the
desired numerical accuracy.
The classical Ewald summation separates the series into long-range smooth parts and short-range singular parts: 
\begin{gather}
\frac{1}{r}=\frac{\erf(\sqrt{\alpha}r)}{r}+\frac{\erfc(\sqrt{\alpha}r)}{r},
\end{gather}
where $\erf(x):=\frac{2}{\sqrt{\pi}}\int_0^x \exp(-u^2)du$ is the error function and $\erfc=1-\erf$. 
Correspondingly,
\begin{gather}
U_1=\frac{1}{2}\sum_{\bm{n}}{}'\sum_{i,j}q_iq_j\frac{\erf(\sqrt{\alpha}|\bm{r}_{ij}+\bm{n}L|)}{|\bm{r}_{ij}+\bm{n}L|}, \\
~~U_2=\frac{1}{2}\sum_{\bm{n}}{}'\sum_{i,j}q_iq_j\frac{\erfc(\sqrt{\alpha}|\bm{r}_{ij}+\bm{n}L|)}{|\bm{r}_{ij}+\bm{n}L|}.
\end{gather}
 The computation of force can be done directly using
\[
\bm{F}_i=-\nabla_{\bm{r}_i}U=-\nabla_{\bm{r}_i}U_1-\nabla_{\bm{r}_i}U_2=:\bm{F}_{i,1}+\bm{F}_{i,2}.
\]
The second part $\bm{F}_{i2}$ corresponds to the short-range forces whose computational cost is relatively low, since, for each particle,  one just needs to
add a finite number of particles in its close neighbour.  We now focus on 
the first part.

The slow decay of $U_1$ in $r$, corresponding to the long-range,  can be dealt with in the Fourier space thanks to its smoothness (see \cite[Chap. 12]{frenkel2001understanding}):
\begin{gather}
U_1=\frac{2\pi}{V}\sum_{\bm{k}\neq 0}\frac{1}{k^2}|\rho(\bm{k})|^2
e^{-k^2/4\alpha}-\sqrt{\frac{\alpha}{\pi}}\sum_{i=1}^N q_i^2,
\end{gather}
where $k=|\bm{k}|$ and $\rho(\bm{k})$ is given by
$\rho(\bm{k}):=\sum_{i=1}^N q_i e^{i\bm{k}\cdot\bm{r}_i}$.
The divergent $\bm{k}=0$ term is usually neglected in simulations to represent that the periodic system is embedded in a conducting medium which is essential for simulating ionic systems.
Then
\begin{gather}\label{eq:force}
\bm{F}_{i,1}=-\sum_{\bm{k}\neq 0}\frac{4\pi q_i \bm{k}}{V k^2}
e^{-k^2/(4\alpha)}\mathrm{Im}(e^{-i\bm{k}\cdot\bm{r}_i}\rho(\bm{k})),
\end{gather}
where we recall $\bm{r}_{ij}=\bm{r}_j-\bm{r}_i$, pointing towards particle $j$ from particle $i$.
Note that the force $\bm{F}_{i,1}$ is bounded for small $\bm{k}$. In fact, $k \ge  2\pi/L$, so $Vk\ge 2\pi L^2$.
Let us consider the factor $e^{-k^2/(4\alpha)}$, and  denote the sum of such factors by
\begin{gather}\label{eq:S}
S:=\sum_{\bm{k}\neq 0}e^{-k^2/(4\alpha)}=H^3-1,
\end{gather}
where
\begin{gather}
H:=\sum_{m\in \bbZ}e^{-\pi^2 m^2/(\alpha L^2)}
=\sqrt{\dfrac{\alpha L^2}{\pi}}\sum\limits_{m=-\infty}^{\infty}e^{-\alpha m^2L^2}
\approx\sqrt{\frac{\alpha L^2}{\pi}}(1+2e^{-\alpha L^2}), \label{psf}
\end{gather}
since often $\alpha L^2 \gg 1$. Hence, $S$ is the sum for all three-dimensional vectors $\bm{k}$ except $0$. 
Then, one can regard the sum as an expectation over the probability distribution
\begin{gather}\label{eq:probexpression}
\mathscr{P}_{\bm{k}}:=S^{-1}e^{-k^2/(4\alpha)},
\end{gather}
which, with $\bm{k}\neq 0$, is a discrete Gaussian distribution and can be sampled efficiently. For example, one can use the Metropolis-Hastings (MH) algorithm (see \cite{hastings1970monte} for details) by choosing proposal samples from the continuous Gaussian
$ \mathcal{N}(0, \alpha L^2/(2\pi^2))$, the normal distribution with mean zero and variance $\alpha L^2/(2\pi^2)$. It should be
emphasized that this sampling can be done {\it offline}, before the 
iteration begins. Once the time evolution starts one just needs to randomly
draw a few ($p$) samples for each time step from this pre-sampled Guassian sequence.

With this observation, the MD simulations can then be done via the random mini-batch approach with this importance sampling strategy. Specifically, one  approximates the force $\bm{F}_{i,1}$ in \eqref{eq:force} by the following random variable:
\begin{gather}\label{eq:rbmapprox}
\bm{F}_{i,1}\approx \bm{F}_{i,1}^*:=-\sum\limits_{\ell=1}^p \dfrac{S}{p}\dfrac{4\pi \bm{k}_\ell q_i}{V k_\ell^2}\mathrm{Im}(e^{-i\bm{k}_\ell\cdot\bm{r}_i}\rho(\bm{k}_\ell)).
\end{gather}
The corresponding algorithm is shown in Algorithm \ref{RBEalg}.
\begin{algorithm}[H]
	\caption{(Random-batch Ewald)}\label{RBEalg}
	\begin{algorithmic}[1]
		\State Choose $\alpha$, $r_c$ and $k_c$ (the cutoffs in real and Fourier spaces respectively), $\Delta t$, and batch size $p$.  Initialize the positions and velocities of charges $\bm{r}^0_i, \bm{v}^0_i$ for $1\le i\le N$.
                 \State Sample sufficient number of $\bm{k}\sim e^{-k^2/(4\alpha)},\, \bm{k}\neq 0$ by the MH procedure to form a set $\mathcal{K}$.
		\For {$n \text{ in } 1: N$}
		\State Integrate Newton's equations \eqref{eq:md1} for time $\Delta t$ with appropriate integration scheme and some appropriate thermostat. The Fourier parts of the Coulomb forces are computed using RBE force \eqref{eq:rbmapprox} with the $p$ frequencies chosen from $\mathcal{K}$ in order.
		\EndFor
	\end{algorithmic}
\end{algorithm}

 Similar to the strategy in the PPPM, one may choose $\alpha$ such that the time cost in real space is cheap and then speed up the computation in the Fourier space. Compared with PPPM, the only difference is that PPPM uses FFT while RBE uses random mini-batch to speed up the computation in the Fourier space. Hence, we make the same choice
\[
\sqrt{\alpha}\sim \frac{N^{1/3}}{L}=\rho_r^{1/3},
\]
which is inverse of the average distance between two numerical particles.
The complexity for the real space part is $\mathcal{O}(N)$. 
 By choosing {\it the same batch} of frequencies for all forces \eqref{eq:rbmapprox} (i.e., using the same
$\bm{k}_{\ell}$, $1\le \ell \le p$ for all $\bm{F}^*_{i,1}, 1\le i\le N$) in the same time step, the complexity per iteration for the frequency part is reduced to $\mathcal{O}(pN)$. This implies that the RBE method has linear complexity per time step if one chooses $p=\mathcal{O}(1)$.

\begin{figure}[ht]	
	\centering
	\includegraphics[width=0.48\textwidth]{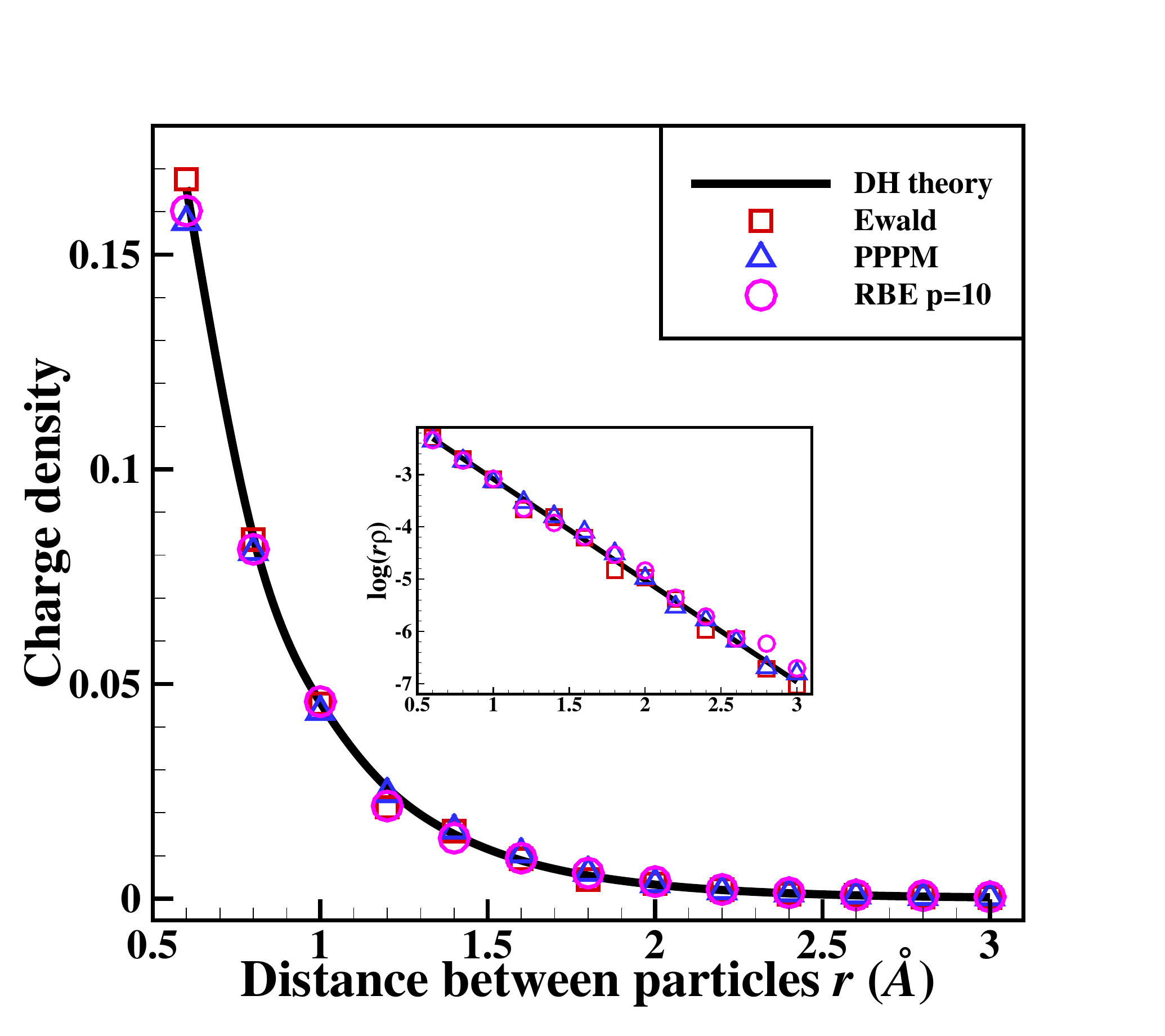}
	\includegraphics[width=0.48\textwidth]{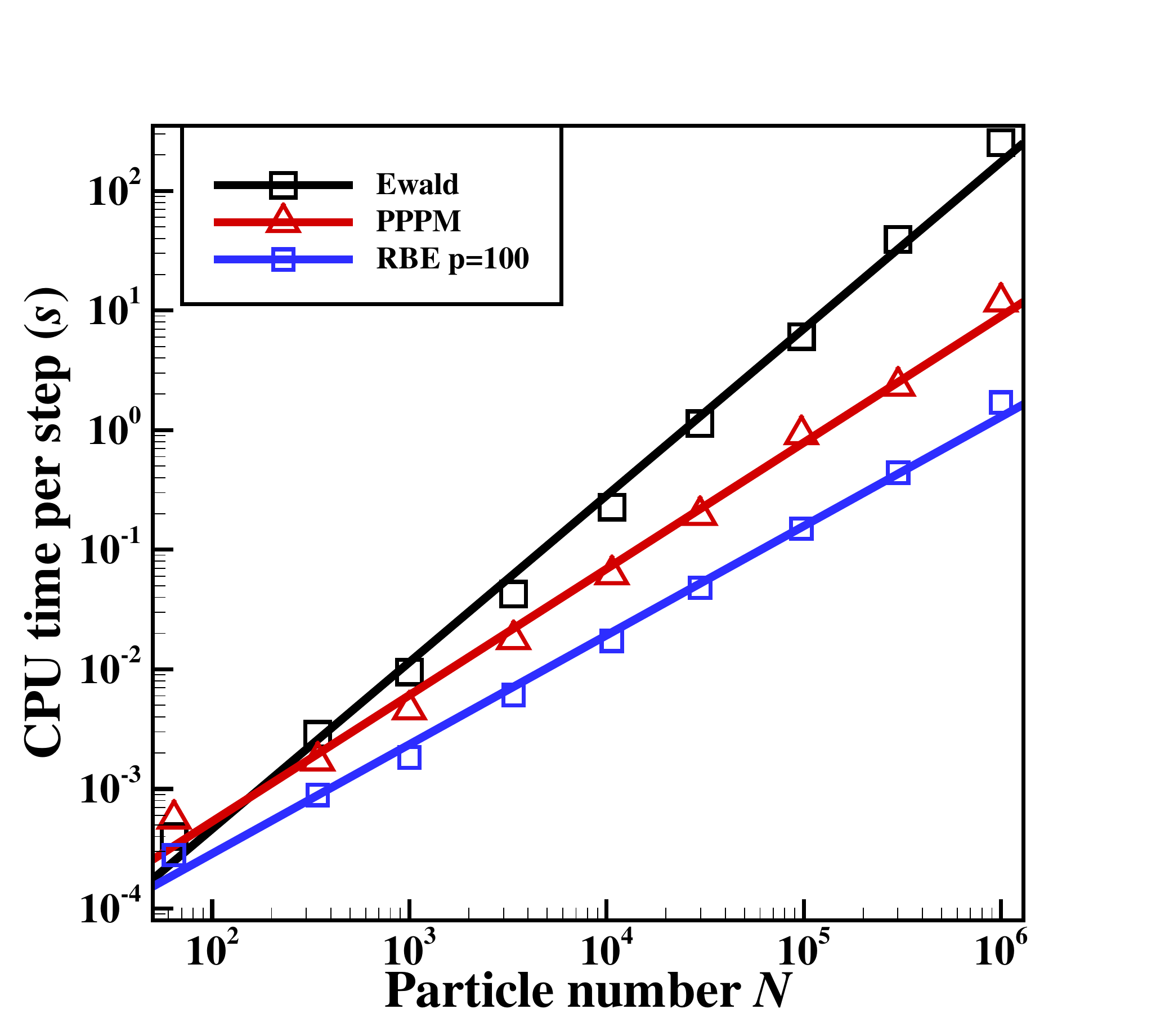}
	\caption{Comparision of the Ewald sum, the PPPM and the RBE methods}
	\label{fig:rbecompare}
\end{figure}

To illustrate the performance of the RBE method, consider an electrolyte with monovalent binary ions (first example in \cite{jinlixuzhao2020rbe}). In the reduced units (\cite[section 3.2]{frenkel2001understanding}), the dielectric constant is taken as $\varepsilon=1/4\pi$ so that the potential of a charge is $\phi(r)=q/r$ and the temperature is $T=\beta^{-1}=1$. Under the Debye-H\"uckel (DH) theory (linearized Poisson-Boltzmann equation), the charge potential outside one ion is given by
\[
-\varepsilon \Delta \phi= 
\begin{cases}
0  & r<a\\
q\rho_{\infty,+}e^{-\beta q\phi}-q\rho_{\infty,-}e^{\beta q\phi}\approx \beta q^2\rho_r \phi, & r>a
\end{cases}
\]
where $\rho_{\infty,+}=\rho_{\infty,-}=N/(2V)$ are the densities of the positive and negative ions at infinity, both being $\rho_r/2$. The parameter $a$ is the effective diameter of the ions, which is related to the setting of the Lennard-Jones potential. In the simulations, $a=0.2$ and the setting of Lennard-Jones potential can be found in \cite{jinlixuzhao2020rbe}.  This approximation gives the net charge density $\rho=-\varepsilon\Delta\phi$
for $r\gg a$,
\[
\ln(r\rho(r))\approx -1.941 r-1.144.
\]
The results in the left panel of Fig.\ref{fig:rbecompare} were obtained by $N=300$ (i.e., $150$ cation and anion particles respectively) in a periodic box with side length $L=10$.
The thermostat was Andersen's thermostat with collision frequency $\nu=3$.  These parameters are chosen such that they give comparable results. Clearly, all the three methods give correct results, agreeing with the curve predicted by the DH theory. Regarding the efficiency, the right figure shows the time consumed for different particle numbers inside the box with the same side length $L=10$. Both the PPPM and RBE methods scale linealy with the particle numbers. However, even for batch size $p=100$, the RBE method consumes much less time. 
The relative accuracies of the potential obtained by RBE against the PPPM are listed in Table \ref{tabl:relaerr}, for different densities $\rho_r=N/L^3$. Clearly, the RBE method has the same level of accuracy compared with the PPPM method for the densities considered.

\begin{table}[ht]
	\centering
	\begin{tabular}{|c|c|c|c|c|}
		\hline
		&$p=10$&$p=20$&$p=50$&$p=100$\\\hline
		$\rho_r=0.1$&$0.15\%$&$0.13\%$&$0.13\%$&$0.08\%$\\\hline
		$\rho_r=0.3$&$0.10\%$&$0.08\%$&$0.04\%$&$0.09\%$\\\hline
		$\rho_r=1$&$0.66\%$&$0.18\%$&$0.11\%$&$0.04\%$\\\hline
		$\rho_r=4$&$7.83\%$&$2.38\%$&$0.71\%$&$0.31\%$\\\hline
	\end{tabular}
	\caption{Relative error of potential energy for the RBE method against PPPM method with different densities and batch sizes.}
	\label{tabl:relaerr}
\end{table}

\begin{figure}[ht]	
	\centering
	\includegraphics[width=0.48\textwidth]{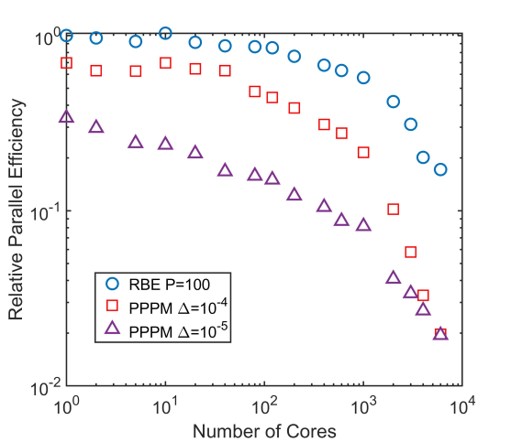}
	\includegraphics[width=0.48\textwidth]{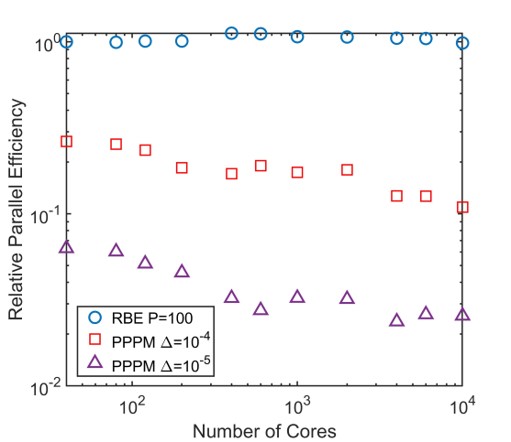}
	\caption{The parallel efficiency of the PPPM and the RBE methods for all-atom simulation of pure water system (Left)  $3\times 10^5$ atoms; (Right) $3\times 10^7$ atoms}
	\label{fig:paraefficiency}
\end{figure}

Next, in Fig. \ref{fig:paraefficiency}, the parallel efficiency of the PPPM and RBE methods from  \cite{liangetalRBE} for the all-atom simulation of pure water systems is shown.  As can be seen, due to the reduction of communications for the particles, the RBE method gains better parallel efficiency. This parallel efficiency is more obvious when the number of particles is larger.  In \cite{liangetalRBE}, the simulation results of pure water system also indicate that the RBE type methods can not only sample from the equilibrium distribution, but also compute accurately the dynamical properties of the pure water systems.

\section{Statistical sampling}

Sampling from a complicated or even unknown probability distribution is crucial in many applications, including numerical integration for statistics of many-body systems \cite{frenkel2001understanding, weinan2019applied}, parameter estimation for Bayesian inference \cite{tarantola2005inverse,box2011bayesian} etc.. The methods that rely on random numbers for sampling and numerical simulations are generally called the Monte Carlo (MC) methods \cite{kalos2009monte,weinan2019applied}. The law of large numbers \cite{durrett2010} validates the usage of empirical measures for approximation of the complicated or unknown probability measure. By the central limit theorem \cite{durrett2010}, the error of the MC methods scales like $\cO(N^{-1/2})$ which is independent of the dimension $d$, hence the MC methods overcome the curse of dimensionality. The Markov Chain Monte Carlo (MCMC) methods \cite{gilks1995markov,gamerman2006markov} are among the most popular MC methods. By constructing Markov chains that have the desired distributions to be the invariant measures, one can obtain samples from the desired distributions by recording the states of the Markov chains.
A typical MCMC algorithm is the Metropolis-Hastings algorithm \cite{metropolis1953,hastings1970monte}. 

Unlike the MCMC, the Stein variational Gradient method (proposed by Liu and Wang in \cite{liu2016stein}) belongs to the class of {\it particle based} variational inference sampling methods (see also \cite{rezende2015variational,dai2016provable}). These methods update particles by solving optimization problems, and each iteration is expected to make progress toward the desired distribution.  As a non-parametric variational inference method, SVGD gives a deterministic way to generate points that approximate the desired probability distribution by solving an ODE particle system, which displays different features from the Monte Carlo methods. 

We describe in this section two sampling methods that use RBM to improve the efficiency. The first method, Random Batch Monte Carlo, is a fast MCMC that costs only $\cO(1)$ per iteration to sample from the Gibbs measures corresponding to many-body particle systems with singular interacting kernels. The second method, RBM-SVGD, is an interesting application of RBM to the Stein variational gradient descent ODE system, which is an interacting particle system.

\subsection{Random Batch Monte Carlo for many body systems}

Suppose that one wants to sample from the $N$-particle Gibbs distribution
\begin{gather}\label{eq:Gibbs}
\pi(\underline{x})\propto \exp\left[-\beta H(\underline{x})\right],
\end{gather}
with $\underline{x}=(x_1,\cdots, x_N)\in \R^{Nd}$ ($x_i\in \R^d$,  and $d\ge 1, d\in \mathbb{N}$), $\beta$ being a positive constant, the $N$-body energy
\begin{gather}
H(\underline{x}):=\sum_{i=1}^N w_i V(x_i)+\sum_{i,j: i<j} w_i w_j \phi(x_i-x_j),
\end{gather}
and $V$ being the external potential assumed to be smooth. Here, $w_i$'s are the weights. In the molecular regime, $w_i$'s are often taken to be $1$, while in the mean field regime \cite{stanley1971, georges1996, lasry2007},
one may have $w\sim N^{-1}$. 

In  \cite{li2020random}, Li et. al. proposed the Random Batch Monte Carlo method, 
which costs $\cO(1)$ per time step for sampling from equilibrium distributions (Gibbs measures) corresponding to particle systems with singular interacting kernels. Similarly to \cite{martin1998novel,hetenyi2002multiple} and the MD methods above, the interacting potential is  decomposed into two parts
\begin{gather}
\phi(x)=\phi_1(x)+\phi_2(x),
\end{gather}
where we suppose that $\phi_1$ has long range but is smooth and bounded, while $\phi_2$ is singular and of short range. The algorithm is based on the following splitting Monte Carlo, which is a special case of the Metropolis-Hastings algorithm:

Suppose there are $N$ particles located at $x_j$ for $j=1,\cdots, N$.  Let us consider the following method for a Markovian jump.

{\it Step 1 ---} Randomly choose a particle $i$.

{\it Step 2 ---}
 Move the particle using $\phi_1$ with overdamped Langevin equation:
 \begin{gather}\label{eq:sde1}
 \begin{split}
 & d\bm{r}_i= -\left(\frac{\nabla V(\bm{r}_i)}{w(N-1)}+\frac{1}{N-1}\sum_{j: j\neq i} \nabla \phi_1(\bm{r}_i-\bm{r}_j)\right) \,dt+\sqrt{\frac{2}{(N-1)w^2 \beta}}\, d\bm{W}_i,\\
 & \bm{r}_i(0)=x_i,
 \end{split}
 \end{gather}
 where $x_j$'s are fixed.
  Evolve this SDE with some time $t>0$ and obtain $\bm{r}_i(t)\rightarrow x_i^*$ as a candidate position of particle $i$ for the new sample.

{\it Step 3 ---}
Use $\phi_2$ to do the Metropolis rejection. Define
 \begin{gather}\label{eq:acceptance}
 \text{acc}(x_i, x_i^*)=\min\left\{1, \exp\Big[-\beta \sum_{j: j\neq i}w^2(\phi_2(x_i^*-x_j)-\phi_2(x_i-x_j))\Big]  \right\}.
 \end{gather}
 With probability $\text{acc}(x_i, x_i^*)$,  accept $x_i^*$ and set
 \begin{gather}
 x_i \leftarrow x_i^*.
 \end{gather}
Otherwise, $x_i$ is unchanged. Then, a new sample $\{x_1,\cdots,x_N\}$ is obtained for the Markov chain.

Note that the overdamped Langevin equation satisfies the detailed balance condition so the above algorithm is a special case of the Metropolis-Hastings algorithm,  thus can correctly sample from the desired Gibbs distribution. 
Due to the short range of $\phi_2$, {\it Step 3} can be done
in $\cO(1)$ operations using some standard data structures such as the cell list \cite[Appendix F]{frenkel2001understanding}. The idea is to use the random mini-batch approach to {\it Step 2}. Hence, one discretizes the SDE with the Euler-Maruyama scheme \cite{kloeden2013numerical,milstein2013stochastic}.
The interaction force is approximated  within the random mini batch idea. This gives the following algorithm.

\begin{algorithm}[H]
\caption{(Random-batch Monte Carlo algorithm)}\label{fastmcmc}
\begin{algorithmic}[1]
\State Split $\phi:=\phi_1+\phi_2$ such that $\phi_1$ is smooth and with long range; $\phi_2$ is with short range.
Generate $N$ initial particles; choose $N_s$ (the total number of samples), $p>1$, $m\ge1$
\For {$n \text{ in } 1: N_s$}
\State Randomly pick an index $i\in \{1, \cdots, N\}$ with uniform probability
       \State $\bm{r}_i \leftarrow x_i$
     \For {$k=1,\cdots, m$}
      \State  Choose $\bm{\xi}_k$, $\bm{z}_k\sim \mathcal{N}(0, I_d)$, $\dt_k>0$ and let,
 \begin{gather*}
 \bm{r}_i\leftarrow \bm{r}_i-\dt_k \left[\frac{\nabla V(\bm{r}_i)}{w(N-1)}+\frac{1}{p-1}\sum_{j\in \bm{\xi}_k} \nabla \phi_1(\bm{r}_i-x_j)\right]+\sqrt{\frac{2\dt_k}{(N-1)w^2 \beta}} \bm{z}_k
 \end{gather*}
      \EndFor
      \State Let $x_i^*\leftarrow \bm{r}_i$. Compute the following using cell list or other data structures:
      \[
      \alpha=\min\left\{1, \exp\Big[-\beta\sum_{j: j\neq i}  w^2(\phi_2(x_i^*-x_j)-\phi_2(x_i-x_j))\Big]  \right\}
      \]
      \State Generate a random number $\zeta$ from uniform distribution on $[0, 1]$. If $\zeta\le \alpha$, set
      \[
      x_i \leftarrow x_i^*
      \]
\EndFor
\end{algorithmic}
\end{algorithm}

It has been proved in \cite{li2020random} that the mini-batch approximation has the an  error control for the transition probability so that the method is correct with some systematic error. The computational cost is $\cO(1)$ for each iteration and the efficiency could be higher since there is no rejection in {\it Step 2}.

We now present a numerical result from \cite{li2020random} to illustrate the efficiency of RBMC.
Consider the Dyson Brownian motion \cite{erdos2017} :
\begin{gather}\label{e:dyson1}
d\lambda_j(t) =-\lambda_j(t)\,dt+\frac{1}{N-1}\sum_{k: k\neq j}\frac{1}{\lambda_j-\lambda_k}dt
+\frac{1}{\sqrt{N-1}} dW_j,~~j=1,\cdots, N,
\end{gather}
where $\{\lambda_j\}$'s represent the eigenvalues of certain random matrices (compared with the original Dyson Brownian motion, $N-1$ instead of $N$ is used in \eqref{e:dyson1}; there is little effect due to the replacement $N\to N-1$). 
In the limit $N\to\infty$,  the distribution obeys the following nonlocal PDE
\begin{gather}\label{eq:dysonlimiteq}
\partial_t\rho(x,t)+\partial_x(\rho(u-x))=0, ~~u(x, t)=\pi(H\rho)(x, t)=\mathrm{p.v.}\int_{\mathbb{R}}\frac{\rho(y,t)}{x-y}\,dy,
\end{gather}
where $H(\cdot)$ is the Hilbert transform on $\mathbb{R}$, $\pi=3.14\cdots$ is the circumference ratio and p.v. represents the Cauchy principal value. From this PDE, one finds that the limiting equation \eqref{eq:dysonlimiteq} has an invariant measure, given by the semicircle law:
\begin{gather}\label{eq:semicircle}
\rho(x)=\frac{1}{\pi}\sqrt{2-x^2}.
\end{gather}

\begin{figure}[ht]	
	\centering
	\includegraphics[width=0.45\textwidth]{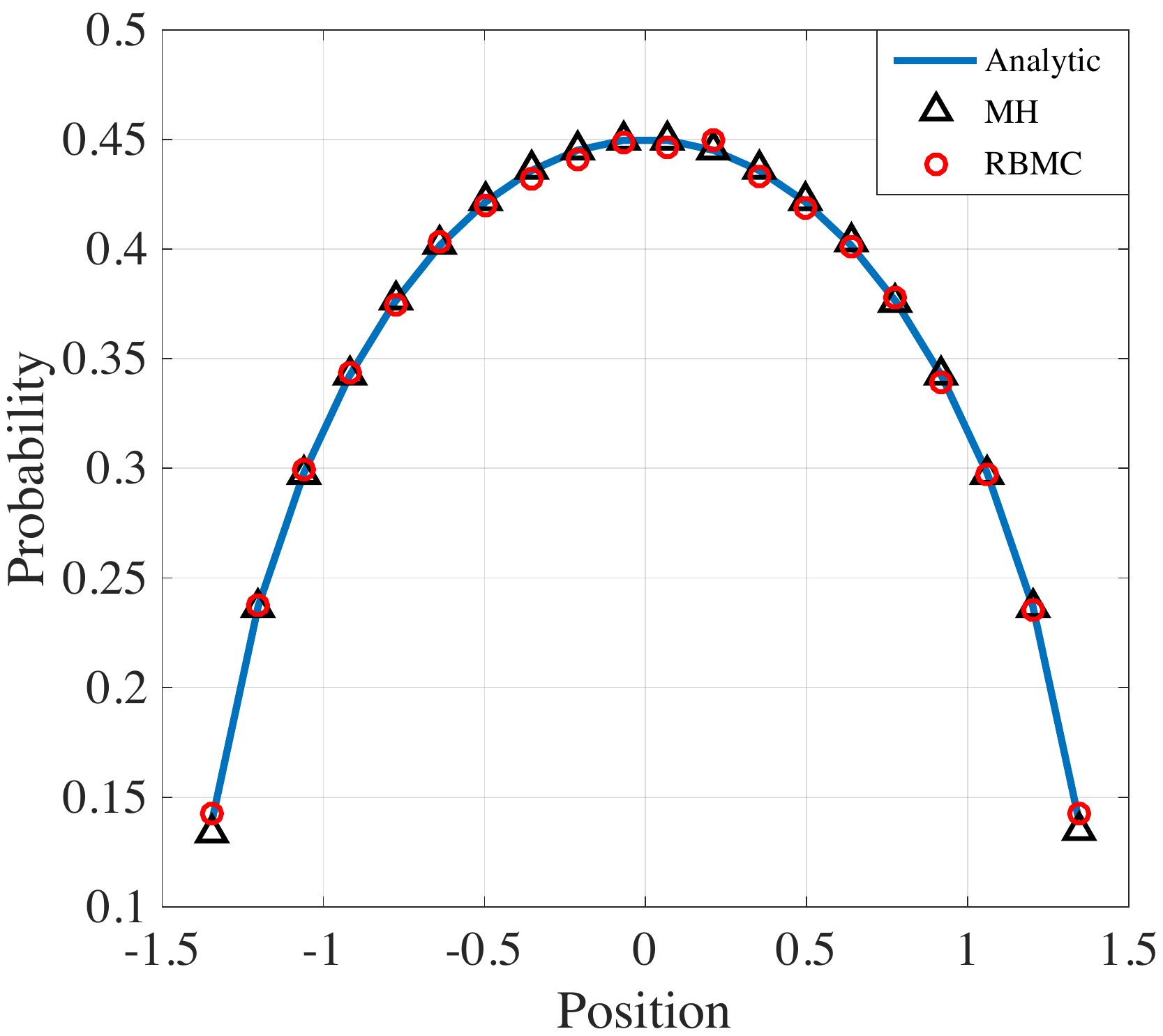}
	\includegraphics[width=0.45\textwidth]{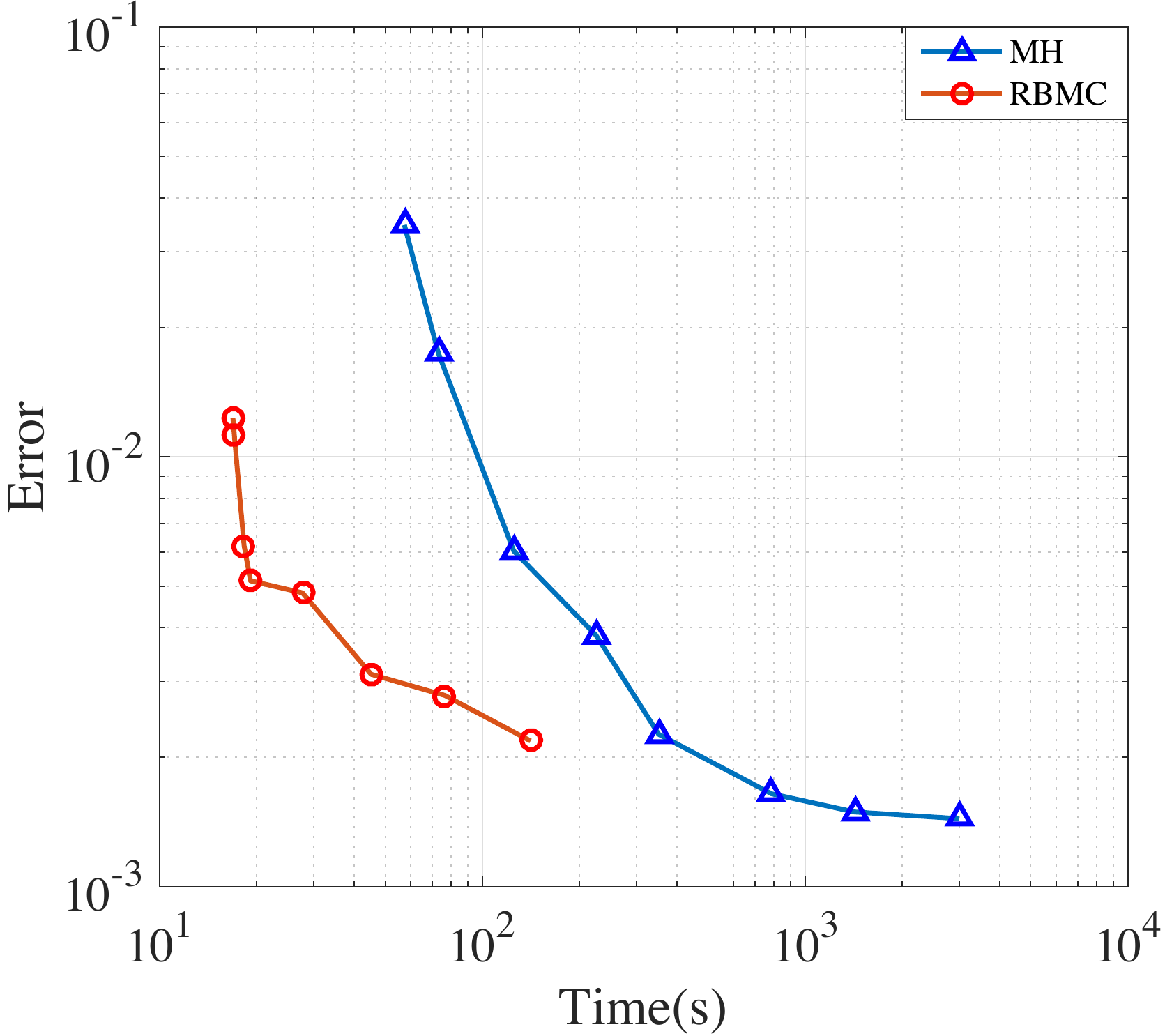}
	\caption{(Left) Empirical densities with 1e7 sampling iterations (1e7$N$ sample points). The blue curve is the analytical curve given by the semicircle law \eqref{eq:semicircle}. (Right) error versus CPU time. }
	\label{fig:dyson}
\end{figure}

Fig.\ref{fig:dyson} shows the sampling results of RBMC and MH methods for empirical measures with particles from the joint distribution
\[
\pi(d\underline{x})\propto \exp\left(-\Big(\frac{N-1}{2}\sum_i x_i^2-\sum_{i<j}\ln|x_i-x_j|\Big)\right),
\]
which is the invariant measure for the interacting particle system \eqref{e:dyson1}. The empirical measure is expected to be close to the semicircle law when $N$ is large enough. In the simulations, the particle number was fixed as $N=500$. In the RBMC, the splitting was done for $\ln r$ at $r=0.01$, and the time step was chosen as $\dt=10^{-4}$. The MH algorithm uses a certain Gaussian proposal for the random movement of a chosen particle. The left panel of Fig. \ref{fig:dyson} shows that both methods yield results that agree with the semicircle law reasonably well. The right panel plots the relative error with respect to the semicircle law versus CPU time.
Clearly, the RBMC method only needs $10\%$ of the time for the MH method to get the error tolerance considered.

\subsection{RBM-SVGD: a stochastic version of stein variational gradient descent}

Suppose that one is  interested in some target probability distribution with density $\pi(x)$ ($x\in\mathbb{R}^d$). In SVGD, one sets $V=-\log \pi$, chooses some symmetric positive definite kernel $\mathcal{K}(x, y)$ and solves the following ODE system for given initial points $\{\bm{r}_i(0)\}_{i=1}^N$ (see \cite{liu2016stein,liu2017stein}):
\begin{gather}\label{eq:discreteODE}
\dot{\bm{r}}_i=
\frac{1}{N}\sum_{j=1}^N \nabla_y \mathcal{K}(\bm{r}_i,\bm{r}_j)
-\frac{1}{N}\sum_{j=1}^N \mathcal{K}(\bm{r}_i,\bm{r}_j)\nabla V(\bm{r}_j),~~i=1,\cdots, N,
\end{gather}
where $N$ is the number of particles for the sampling purpose. The subindex ``$y$'' in $\nabla_y$ means that the gradient is taken with respect to the second variable in $\mathcal{K}(\cdot, \cdot)$; i.e. $\nabla_y \mathcal{K}(\bm{r}_i,\bm{r}_j):=\nabla_y\mathcal{K}(x, y)|_{(x, y)=(\bm{r}_i, \bm{r}_j)}$. When $t$ is large enough, the empirical measures constructed using $\{\bm{r}_i(t)\}_{i=1}^N$ is expected to be close to $\pi$, i.e.
\[
\frac{1}{N}\sum_{i=1}^N \delta(x-\bm{r}_i(t)) \approx \pi(x)\,dx,~~t\gg 1.
\]
SVGD provides consistent estimation for generic distributions as Monte Carlo methods do, but it seems to be more efficient than some Monte Carlo methods in practice level for approximating the desired measure, when the number of particles is small \cite{liu2016stein,detommaso2018stein}. Interestingly, it reduces to the maximum a posterior (MAP) method when $N=1$ \cite{liu2016stein}. 

The ODE system \eqref{eq:discreteODE} clearly is an interacting particle system but now the interaction kernel is no longer translation invariant and is not symmetric. The kernel can even grow as $|\bm{r}_i-\bm{r}_j|\to\infty$. Clearly, for such systems, RBM is applicable. Applying the RBM to this special kernel and using any suitable ODE solvers, one gets a class of sampling algorithms, which is called RBM-SVGD  in \cite{li2020svgd}.  The discrete algorithm (with possible variant step size) is shown in Algorithm~\ref{ssvgd}. Clearly, the complexity is $\cO(pN)$ for each iteration.

\begin{algorithm}[H]
    \caption{RBM-SVGD}\label{ssvgd}
    \begin{algorithmic}[1]
        \For {$k \text{ in } 0: N_T-1$}   
            \State Divide $\{1, 2, \ldots, pn\}$ into $n$ batches randomly.
            \For {each batch  $\mathcal{C}_q$} 
                \State For all $i\in \mathcal{C}_q$,
                \[
                    \bm{r}_i^{(k+1)}\leftarrow \bm{r}_i^{(k)}
                    + \frac{1}{N} \Big( \nabla_y
                    \mathcal{K}(\bm{r}_i^{(k)}, \bm{r}_i^{(k)})
                    - \mathcal{K}(\bm{r}_i^{(k)}, \bm{r}_i^{(k)}) \nabla V(\bm{r}_i^{(k)})\Big)\eta_{k}
                    +\Phi_{k,i}\eta_k,
                \]
                where
                \begin{gather}
                    \Phi_{k,i} = \frac{N-1}{N(p-1)}
                    \sum_{j\in \mathcal{C}_q,j\neq i}
                    \left(\nabla_y \mathcal{K}(\bm{r}_i^{(k)}, \bm{r}_j^{(k)})
                    -\mathcal{K}(\bm{r}_i^{(k)}, \bm{r}_j^{(k)}) \nabla V(\bm{r}_j^{(k)})\right).
                \end{gather}
            \EndFor
        \EndFor
    \end{algorithmic}
\end{algorithm}
Here, $N_T$ is the number of iterations and $\{\eta_k\}$ is the sequence of time steps, which play the same role as learning rate in SGD \cite{bottou1998online,bubeck2015convex}. For some applications,  one may simply set $\eta_k=\eta\ll 1$ to be a constant and get relatively good results. However, in many high dimensional problems, choosing $\eta_k$ to be constant may yield divergent sequences  \cite{robbins1951stochastic}. One may decrease $\eta_k$ to obtain convergent data sequences. For example, one may simply choose $\eta_k=1/k$ as in SGD. Another frequently used strategy is the AdaGrad approach \cite{duchi2011adaptive,ward2019adagrad}.


We recall the gradient flow under the so-called  ``Stein metric" in the space of probability measures \cite{liu2017stein,gao2020note}: 
\begin{align}\label{eq:odeflow} 
\partial_t\rho=\nabla\cdot\left(\rho \mathcal{K}*( \rho \nabla\frac{\delta E}{\delta\rho}) \right),
\end{align}
where $\mathcal{K}*g=\int\mathcal{K}(x, y)g(y)\,dy$. Consider taking the energy functional as the Kullback-Leibler (KL) divergence between $\rho$ and the target distribution $\pi$, where KL divergence is also known as the relative entropy defined by
\begin{gather}
    \kl(\mu || \nu)=\mathbb{E}_{Y\sim
    \mu}\log\left(\frac{d\mu}{d\nu}(Y)\right).
\end{gather}
Here $\frac{d\mu}{d\nu}$ is the well-known Radon-Nikodym derivative.
Then, equation (\ref{eq:odeflow}) becomes
\begin{gather}\label{eq:svgdmeanfield}
    \partial_t \rho=\nabla\cdot(\rho \mathcal{K}*(\rho\nabla V+\nabla \rho)).
\end{gather}
It is easy to see that $\pi\propto \exp(-V)$ is invariant under this PDE.
See   \cite{liu2017stein,lu2018scaling} for some relevant studies.

The above theory encounters difficulty for empirical measures because the KL divergence is simply infinity.  One benefit of the of the ``Stein metric" is that the gradient may be moved from $\nabla\rho$ onto the kernel $\mathcal{K}(x,y)$ so that the flow \eqref{eq:odeflow} becomes  \eqref{eq:discreteODE}, which is then well-defined. In fact, if $\{\bm{r}_i\}$ solves
the ODE system \eqref{eq:discreteODE}, then the corresponding empirical measure is a measure solution to \eqref{eq:svgdmeanfield} (see \cite[Proposition 2.5]{lu2018scaling}).  Hence, one may reasonably expect that \eqref{eq:discreteODE} will give approximation for the desired distribution $\pi$.

\begin{figure}[ht]	
	\centering
	\includegraphics[width=0.6\textwidth]{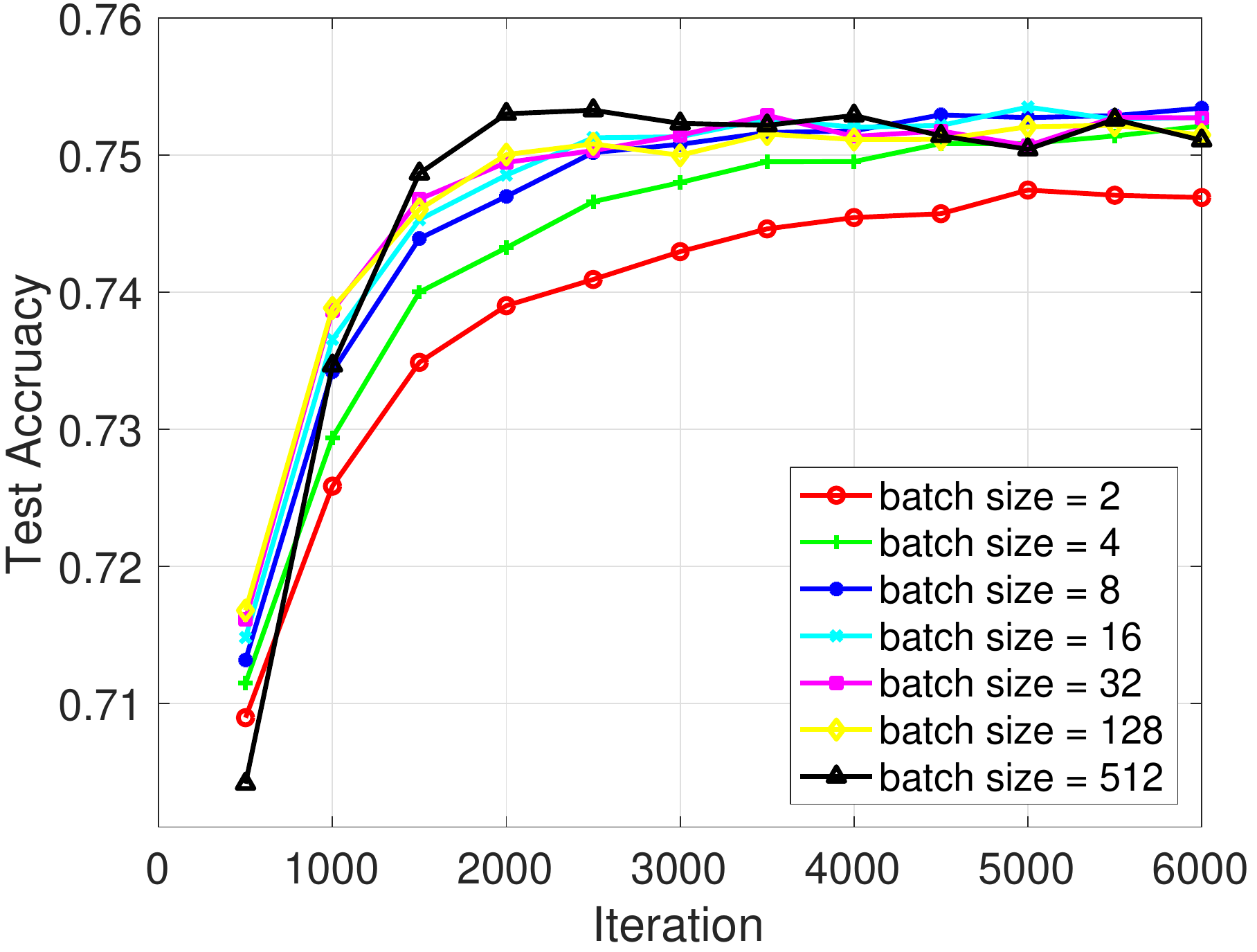}
	\caption{Test accuracy on the Covertype dataset}
	\label{fig:bayes}
\end{figure}

For numerical illustration, we take an example from  \cite{li2020svgd}.  Consider the logistic regression for binary classification
on the Covertype dataset,  with 581012 data points and 54 features~\cite{gershman2012nonparametric}.
The inference is applied on posterior $p(x|D)$ with the parameter $x=[w,\log\alpha]$ being of dimension 55. 
 Here, $D$ is $80\%$ of the data and the remaining data were used for test.
 Figure \ref{fig:bayes} shows the performance of SVGD  and RBM-SVGD
with $N=512$ particles and kernel $\mathcal{K}(x,y)=k(x-y)$ for a Gaussian kernel $k(\cdot)$.
Clearly, RBM-SVGD gives comparable results with SVGD, both results being as good as some traditional methods.

Table \ref{tab:bayselr-runtime} shows the CPU time and speedup of RBM-SVGD. Clearly, for comparable results,
RBM-SVGD is more efficient.

\begin{table}[htp]
    \centering
    \caption{Average runtime of $6000$ iterations}
    \label{tab:bayselr-runtime}
    \begin{tabular}{cccccccc}
        \toprule
        & \multicolumn{6}{c}{RBM-SVGD} & SVGD \\
        \cmidrule{2-7}
        p & 2 & 4 & 8 & 16 & 32 & 128 & 512\\
        \midrule
        Runtime(s) & 8.59 & 11.24 & 16.28 & 26.15 & 21.66 & 19.42 &
        47.01 \\
        Speedup & 5.5x & 4.2x & 2.9x & 1.8x & 2.2x & 2.4x & \\
        \bottomrule
    \end{tabular}
\end{table}

\section{Agent-based models for collective dynamics}

Collective behaviors of self-propelled particles (agents) are ubiquitous in nature, for example, synchronous flashing of fireflies and pacemaker cells, swarming of fish, flocking of birds  and herding of sheep. We refer to \cite{A-B-F, cucker2007emergent, motsch2014, T-T, Wi1} for survey articles and related literature.

While the RBM was introduced as an efficient algorithm for interacting particle systems, one can also view it as a (random) model of the underlying problem, which takes into account only
a small number of interactions randomly at discrete time steps. 
Two natural questions arise with such models: (a) How accurate
are these ``new" random models compared to the original, full batch
models? (b) Do these random models still capture the main features
of the original model, such as the collective or long time behavior, and under what conditions? Here we review some recent  results that address these issues for two representative problems, the Cucker-Smale model for flocking and the consensus model.

\subsection{The Cucker-Smale model}

We begin with the Cucker-Smale (CS) model \cite{cucker2007emergent}:
\begin{equation}\label{eq:cs}
\begin{cases}
\displaystyle \frac{d\bm{x}_i}{dt} = \bm{v}_i, \quad t > 0, \quad i=1,\ldots,N, \vspace{0.2em}\\
\displaystyle \frac{d\bm{v}_i}{dt} = \frac{\kappa}{N-1}\sum_{j:j\neq i} \psi(|\bm{x}_j-\bm{x}_i|)(\bm{v}_j-\bm{v}_i),\\
\displaystyle (\bm{x}_i(0), \bm{v}_i(0)) =(x_i^{in}, v_i^{in}), 
\end{cases}
\end{equation}
 where $\bm{x}_i$ and $\bm{v}_i$ are the position and velocity of the $i$-th CS particle, respectively, $\kappa$ is the non-negative coupling strength and $\psi$, the communication weight measuring mutual interactions between agents,  is   positive, bounded, and Lipschitz continuous and satisfies the monotonicity conditions: 
\begin{align}
\begin{aligned} \label{A-2}
& 0 \leq \psi(r) \leq \psi_M, \quad \forall~r \geq 0, \quad \| \psi \|_{\mbox{Lip}} < \infty, \\
& (\psi(r_1) - \psi(r_2)) (r_1 - r_2)  \leq 0, \quad r_1,~r_2 \in \R_+.
\end{aligned}
\end{align}
Here   $\psi_M > 0$ is a constant.  
The emergent dynamics of \eqref{eq:cs}, flocking, in which all particles will eventually stay in a bounded domain with the same velocity,
has been extensively studied in literature \cite{hasimple2009, H-Tad}.

Consider the RBM-approximation for \eqref{eq:cs}: 
\begin{equation}
\begin{cases}
\displaystyle \frac{d\tilde{\bm{x}}_i}{dt}  = \tilde{\bm{v}}_i, \quad  t \in [t_{m-1},t_{m}),~ m=1,2,\cdots,  \\
\displaystyle \frac{d \tilde{\bm{v}}_i}{dt} = \frac{\kappa}{p-1}\sum_{j \in \cC_i^{(k)}, j\not=i} \psi(|\tilde{\bm{x}}_j-\tilde{\bm{x}}_i|)(\tilde{\bm{v}}_j-\tilde{\bm{v}}_i), \\
\displaystyle (\tilde{\bm{x}}_i(0), \tilde{\bm{v}}_i(0)) = (x_i^{in}, v_i^{in}), \quad i=1,\ldots,N.
\end{cases}
\end{equation}
 
 Assume that $\psi$ is long-ranged: 
\begin{equation}\label{A-3-1}
1/\psi(r) = {\mathcal O}(r^\beta) \quad\text{as }~ r \to \infty \quad \mbox{for some $\beta \in [0, 1)$}.
\end{equation}
For example, one can take
\[
\psi(r) = \frac{1}{(1+r^2)^{\beta/2}}, \quad \beta \in [0, 1).
\]
Then \cite{HJKK-C-S-RBM} establishes the following  emergence of a global flocking: there exist positive constants $\tilde{x}_{\infty}$ and $C$ such that 
\begin{align}
\begin{aligned} \label{A-4}
& \sup_{0 \leq t < \infty}{\mathbb E} \Big( \frac{1}{N^2}\sum_{i,j=1}^N|\tilde{\bm{x}}_i - \tilde{\bm{x}}_j|^2 \Big)  < \tilde{x}_{\infty}\quad\text{and}\quad\\
&{\mathbb E} \Big( \frac{1}{N^2} \sum_{i,j=1}^N|\tilde{\bm{v}}_i - \tilde{\bm{v}}_j|^2 \Big) \leq C \exp\Big[  -\frac{C(p-1)}{(N-1)(1+\dt)} t(1+t)^{-\beta} \Big],
\end{aligned}
\end{align}
where $C$ depends only on $\psi$, $\beta$, $\kappa$ and the  initial data.

Furthermore, the following  uniform-in-time error estimate was also proved: when $\psi$ has a positive lower bound $\psi_0$,
\begin{equation}\label{A-5-0}
\begin{aligned}
\psi(r) \geq \psi_0 \quad\text{for }~ r \geq 0,
\end{aligned}
\end{equation}
then 
\begin{align}
\begin{aligned} \label{A-5}
{\mathbb E} \Big( \frac{1}{N}\sum_{i=1}^N  |\tilde{\bm{v}}_i(t) - \bm{v}_i(t)|^2 \Big) &\leq C\dt\left(\frac{1}{p-1}-\frac{1}{N-1}\right)+C\dt^2 \\
& \hspace{2cm} + C (1+\dt)\exp(-\kappa \psi_0 t),
\end{aligned}
\end{align}
where the dependency of the constant $C$ is the same as in \eqref{A-4}. 

Note that the positive lower bound assumption \eqref{A-5-0} corresponds to the case of $\beta = 0$ in the long-ranged communication \eqref{A-3-1}. However, the third time-decaying term in the right-hand side of \eqref{A-5} is independent of $p$ and $N$. 

\subsection{Consensus models}

Let $\bm{q}_i \in \R^d$, $1\le i\le N$ be a collection of agents that  seek for a consensus, governed by the Cauchy problem:
\begin{equation}
\begin{cases} \label{A-0}
\displaystyle \frac{d\bm{q}_i}{dt} =\nu_i+\frac{\kappa}{N-1}\sum_{j \not = i} a_{ij} \Gamma(\bm{q}_j-\bm{q}_i),\quad t > 0, \\
\displaystyle \bm{q}_i(0) = q_i^{in}, \quad  i = 1, \cdots, N,
\end{cases}
\end{equation}
where $\kappa$ is a non-negative coupling strength and $\nu_i$ is the intrinsic velocity of the $i$-th agent. 
Here $\Gamma$ is an interaction function satisfying the following properties: there exists $C_1 > 0$ such that
\begin{equation} 
\Gamma \in {\mathcal C}^2(B_{C_1}(0)), \quad \Gamma(-q) = -\Gamma(q), \qquad \forall~q \in \overline{B_{C_1}(0))}. 
\end{equation}
Here    $B_r(x)$ is the open ball with radius $r$ and center $x$. We  assume, without loss of generality, that the total sum is zero:
\[ \sum_{i=1}^{N} \nu_i = 0, \]
and the adjacency matrix $(a_{ij})_{i,j=1}^N$ represents the network structure for interactions between agents satisfying symmetry and non-negative conditions:
\begin{equation*} \label{A-0-0}
a_{ij} = a_{ji} \geq 0, \quad 1 \leq i, j \leq N.
\end{equation*}

Note that the first term on the R.H.S. of \eqref{A-0} induces the "{\it dispersion effect}" due to the heterogeneity of $\nu_i$.
 The second term in the R.H.S. of \eqref{A-0}, modeled by the convolution type consensus force, generates "{\it concentration effect}",  The overall dynamics of \eqref{A-0} is determined by the competitions between dispersion and concentration.

Below we present the study on RBM to this problem in \cite{K-H-J-K3}. Conisder the RBM-approximation where the interaction term is approximated by the random mini-batch at each time step. Then the relative state $|\tilde{\bm{q}}_i - \tilde{\bm{q}}_j|$ for RBM aproximation can be {\it unbounded} even if the original relative state $|\bm{q}_i - \bm{q}_j|$ is uniformly bounded. Thus to balance dispersion and interaction in the RBM, one also needs to apply the RBM in the dispersion part as well.   A sufficient framework leading to the uniform boundedness of relative states is to introduce suitable decomposition of the dispersion term $\nu_i$ as a sum of $N$-dispersion terms $\bar\nu_{ij}$:
\begin{equation}\label{A-2-1}
\bar\nu_{ij} = -\bar\nu_{ji}, \quad  \nu_i = \frac{\kappa}{N-1}\sum_{j=1}^N\bar\nu_{ij},\quad i,j=1,\dots,N.
\end{equation}
Then, the original Cauchy problem \eqref{A-0} is equivalent to the following problem:
\begin{equation}
\begin{cases} \label{A-2-2}
\displaystyle \frac{d\bm{q}_i}{dt} = \frac{\kappa}{N-1}\sum_{j \not = i} \Big( \nu_{ij} + a_{ij} \Gamma(\bm{q}_j-\bm{q}_i) \Big),\quad t > 0, \\
\displaystyle \bm{q}_i(0) = q_i^{in}, \quad  i = 1, \cdots, N,
\end{cases}
\end{equation}
and the RBM  samples dispersions and interactions proportionally,
\begin{equation}\label{A-3}
\begin{cases}
\displaystyle \frac{d\tilde{\bm{q}}_i}{dt} =\frac{\kappa}{p-1}\sum_{j \in \cC_i^{(k)}, j \not= i} \big(\bar\nu_{ij}+a_{ij} \Gamma(\tilde{\bm{q}}_j-\tilde{\bm{q}}_i)\big), \quad t \in (t_k,t_{k+1}), \\
\displaystyle \tilde{\bm{q}}_i(0) = q_i^{in},\quad i=1,\ldots,N,~~k = 0, 1, 2, \dots.
\end{cases}
\end{equation}

We  first state the main result for the {\it one-dimensional} case. Assume that the coupling function $\Gamma$ is strongly dissipative in the sense that 
\[  (\Gamma({q_1})-\Gamma(q))\cdot({q_1}-q) \approx |{q_1} -q|^2, \quad \forall~q, {q_1} \in [-C_1, C_1], \]
and also the full system \eqref{A-0} has an equilibrium $\Phi = (\phi_1, \cdots, \phi_N) \in (-C_1,C_1)^N$ with initial data  sufficiently close to $\Phi$.
 The main result is the following uniform error estimate, under the condition that the underlying network topology is connected strongly enough:
\[
\sup_{0 \leq t < \infty} \Big[ \frac{1}{N}\sum_{i=1}^N\mathbb E |\tilde{\bm{q}}_i(t) - \bm{q}_i(t)|^2 \Big] 
\lesssim \Big[ \dt\left(\frac{1}{p-1}-\frac{1}{N-1}\right) + \dt^2 \Big ].  
\]
For the multi-dimensional setting with $\bm{q}_i \in \R^d$, the same error analysis can be obtained under one more extra assumption, which guarantees that the states $Q:=(\bm{q}_1, \cdots, \bm{q}_N)$ and $\tilde{Q}:=(\tilde{\bm{q}}_1, \cdots, \tilde{\bm{q}}_N)$ are confined in the symmetric interval. \vspace{1em}

Now, we give two main results on the emergent dynamics of \eqref{A-3} proved in \cite{HJKK3}. Introduce two functionals for $ \tilde{Q}=(\tilde{\bm{q}}_1, \cdots, \tilde{\bm{q}}_N)$:
\begin{equation*} \label{A-6-1}
{\mathcal M}_2(\tilde{Q}) := \frac{1}{N} \sum_{j=1}^N | \tilde{\bm{q}}_j |^2, \quad  {\mathcal D}(\tilde{Q}) := \max_{1 \leq i, j \leq N} |\tilde{\bm{q}}_i - \tilde{\bm{q}}_j |.
 \end{equation*}
The first main result is concerned with the exponential decay of the second moment of ${\hat q}_i^R$: there exists a positive constant $\Lambda_1 = \Lambda_1(N, P, \tau, \kappa, L_1)$ satisfying
\begin{equation*} \label{A-8}
{\mathbb E} \Big({\mathcal M}_2( \tilde{Q}(t))\Big)\leq e^{-\Lambda_1 t} {\mathbb E} \Big({\mathcal M}_2(\tilde{Q}(0)) \Big), \quad t \geq 0.
\end{equation*}

The second main result deals with almost sure (a.s.) convergence of $\tilde{Q}$: there exists a positive constant $\Lambda_2 = \Lambda_2(N,P,\tau,\kappa,L_1,L_2)$ such that 
\[  {\mathcal D}(\tilde{Q}(t)) \leq   {\mathcal D}(\tilde{Q}(0)) Ce^{-\Lambda_2 t}, \quad t \geq 0. \]

We remark that although the exponential decay rates in above results depend on $N$, numerical results in \cite{HJKK3} show that the decay rates are in not sharp, and they are independent of $N$.

\section{Quantum dynamics}

In this section, we have a review of  the applications of RBM to interacting particles in the quantum regime. 
In particular, we first present and comment on the convergence results of RBM applied to the $N$-body Schr\"odinger
equation in \cite{golse2019random}, and then have a review of the application of RBM to quantum Monte Carlo (QMC) methods in \cite{J-L-QMC}.

\subsection{A theoretical result on the $N$-body Schr\"odinger equation}

The first principle computation is based on solving for complex-valued  wave function $\Psi_N\equiv\Psi_N(t,x_1,\ldots,x_N)\in\mathbb{C}$ of the $N$-body Schr\"odinger equation
\be\lb{NBodySchro}
i\hb\d_t\Psi_N(t,x_1,\ldots,x_N)=\cH_N\Psi_N(t,x_1,\ldots,x_N)\,,\quad\Psi_N\rstr_{t=0}=\Psi_N^{in}
\ee
where $t\ge 0$ is the time while $x_m\in\R^d (1\le m\le N)$ is the position of the $m$th particle, 
$\cH_N$ is the quantum Hamiltonian for $N$ identical particles with unit mass:
\be\lb{NBodyHam}
\cH_N:=\sum_{m=1}^N-\tfrac12\hbar^2\Delta_{x_m}+\tfrac1{N-1}\sum_{1\le \ell<n\le N}V(x_\ell-x_n)\,,
\ee
while $\hbar$  is the reduced Planck constant. The $N$ particles in this system interact via a binary (real-valued) potential $V$ assumed to be even, bounded and sufficiently regular (at least of class $C^{1,1}$ on $\R^d$). 
The coupling constant $\tfrac1{N-1}$ is chosen in order to balance the summations in the kinetic energy (involving $N$ terms) and in the potential energy (involving $\tfrac12N(N-1)$ terms). 

When solving  \eqref{NBodySchro}, the computation
is exceedingly expensive due to the 
smallness of $\hb$, which requires small time steps $\Delta t$ and small mesh sizes of order $\hb$ for the convergence of the numerical scheme, due to the oscillation in the wave function $\Psi_N$ with frequency of order $1/\hb$  
(see \cite{B-J-M, J-M-S}). On top of this, any numerical scheme for \eqref{NBodySchro} requires computing, at each time step, the sum of the interaction potential for each particle pair in the $N$-particle system, which needs  $\cO(N^2)$ operations. The RBM described below 
 reduces  the computational cost to $\cO(N)$ per time-step.

Below we follow the presentation of \cite{golse2019random}. Assume for simplicity that $N\ge 2$ is an even integer. Let $\si_1,\si_2,\ldots,\si_j,\ldots$ be a mutually independent and uniformly distributed random sequence of  permutations. Each permutation 
$\si\in\fS_N$ defines a partition of $\{1,\ldots,N\}$ into $N/2$ batches of two indices:
$$
\{1,\ldots,N\}=\bigcup_{k=1}^{N/2}\{\si(2k-1),\si(2k)\}\,.
$$
Set
\begin{equation}
\bT_t(\ell,n):=\left\{\ba{}&1\quad&&\text{ if }\{\ell,n\}=\left\{\si_{[\frac{t}{\Delta t}]+1}(2k\!-\!1),\si_{[\frac{t}{\Delta t}]+1}(2k)\right\}\text{ for some }k\le \tfrac{N}2\,,\\ 
&0&&\text{ otherwise,}\ea\right.
\end{equation}
and consider the time-dependent random batch Hamiltonian
\be\lb{RBHam}
\wtilde{\cH}_N(t):=\sum_{m=1}^N-\tfrac12\hbar^2\Delta_{x_m}+\sum_{1\le \ell<n\le N}\bT_t(\ell,n)V(x_\ell-x_n)\,.
\ee

The  RBM then solves the random batch Schr\"odinger equation
\be\lb{RBSchro}
i\hb\d_t\wtilde\Psi_N(t,x_1,\ldots,x_N)=\wtilde{\cH}_N(t)\wtilde\Psi_N(t,x_1,\ldots,x_N)\,,\quad\wtilde\Psi_N\rstr_{t=0}=\wtilde\Psi_N^{in}\,.
\ee
 Clearly,  for each time step the cost of computing the interaction potential is reduced from 
$\cO(N^2)$ to $\cO(N)$.

As we have seen, RBM is known to converge in the case of classical dynamics. It is therefore natural to seek an error estimate for the quantum RBM method. The major difficulty here is to obtain an error
estimate that is {\it independent} of $\hbar$ and $N$.

\subsubsection{Mathematical Setting and Main Result}

It will be more convenient to carry out the analysis on the corresponding von Neumann equations
\be\lb{NBodyvN}
i\hb\d_tR_N(t)=\cH_NR_N(t)-R_N(t)\cH_N=:[\cH_N,R_N(t)]\,,\quad R_N(0)=R_N^{in}\,.
\ee
Here we denote $\fH:=L^2(\R^d;\mathbb{C})$ and $\fH_N=\fH^{\otimes N}\simeq L^2((\R^d)^N;\mathbb{C})$ for each $N\ge 2$. The algebra of bounded operators on $\fH$ is denoted by $\cL(\fH)$, while $\cL^1(\fH)\subset\cL(\fH)$
and $\cL^2(\fH)$ are respectively the two-sided ideals of trace-class and Hilbert-Schmidt operators on $\fH$. The operator norm of $A\in\cL(\fH)$ is denoted $\|A\|$. A density operator on $\fH$ is a trace-class operator $R$ on 
$\fH$ such that
$$
R=R^*\ge 0\quad\text{ and }\quad\Tr_\fH(R)=1\,.
$$
The set of density operators on a separable Hilbert space $H$ is
henceforth denoted $\cD(H)$.

The random batch von Neumann equation is 
\be\lb{RBvN}
i\hb\d_t\wtilde R_N(t)=[\wtilde{\cH}(t),\wtilde R_N(t)]\,,\quad \quad\wtilde R_N(0)=R_N^{in}\,.
\ee

In order to find an error estimate for the RBM that is independent of the particle number $N$, one first needs to define in terms of $R_N(t)$ and $\wtilde R_N(t)$ {\it quantities of interest} to be compared that are 
{\it independent of $N$}.  A common practice when considering large systems of identical particles is to study the reduced density operators, which unfortunatey does not lead to $N$-independent error estimates \cite{golse2019random}. 
Assume that $R_N^{in}$ has an integral kernel $r^{in}\equiv r^{in}(x_1,\ldots,x_N;y_1,\ldots,y_N)$ satisfying the symmetry
\be\lb{RinSym}
r^{in}(x_1,\ldots,x_N;y_1,\ldots,y_N)=r^{in}(x_{\si(1)},\ldots,x_{\si(N)};y_{\si(1)},\ldots,y_{\si(N)})
\ee
for  each permutation $\si\in\fS_N$. Then, for each $t\ge 0$, the $N$-body density operator $R_N(t)$ solution of \eqref{NBodyvN} satisfies the same symmetry, i.e. its
 integral kernel of the form $r(t;x_1,\ldots,x_N;y_1,\ldots,y_N)$ also satisfies
\be\lb{RtSym}
r(t;x_1,\ldots,x_N;y_1,\ldots,y_N)=r(t;x_{\si(1)},\ldots,x_{\si(N)};y_{\si(1)},\ldots,y_{\si(N)})
\ee
for each permutation $\si\in\fS_N$. The $1$-particle reduced density operator of $R_N(t)\in\cD(\fH_N)$ is $R_{N,\indc}(t)\in\cD(\fH)$ defined by the integral 
kernel
\be\lb{R1t}
r_\indc(t,x,y):=\int_{(\R^d)^{N-1}}r(t;x,z_2,\ldots,z_N;y,z_2,\ldots,z_N)dz_2\ldots dz_N\,.
\ee

Even if $R_N^{in}$ satisfies the symmetry \eqref{RinSym}, in general $\wtilde R_N(t)$ does not satisfy the symmetry analogous to \eqref{RtSym} for $t>0$ (with $r$ replaced with $\wtilde r$, an integral kernel for 
$\tilde{R}_N(t)$) because the random batch potential
$$
\sum_{1\le \ell<n\le N}\bT_t(\ell,n)V(x_\ell-x_n)
$$
is not invariant under permutations of the particle labels.
For that reason, the $1$-particle reduced density operator of $\wtilde R_N(t)$ one needs is $\wtilde R_{N,\indc}(t)\in\cD(\fH)$ defined for all $t>0$ by the integral kernel
\be\lb{wtR1t}
\wtilde r_\indc(t,x,y):=\frac1N\sum_{j=1}^N\int_{(\R^d)^{N-1}}\wtilde r(t;Z_{j,N}[x],Z_{j,N}[y])d\hat Z_{j,N}\,,
\ee
with the notation
$$
Z_{j,N}[x]:=z_1,\ldots,z_{j-1},x,z_{j+1}\ldots,z_N\,,\quad d\hat Z_{j,N}=dz_1\ldots dz_{j-1}dz_{j+1}\ldots dz_N\,.
$$
(Obviously \eqref{wtR1t} holds with $r_\indc$ and $r$ in the place of $\wtilde r_\indc$ and $\wtilde r$ respectively because of the symmetry \eqref{RtSym}.)

\smallskip
We also need to introduce the Wigner functions of the 
density operators $R_N(t)$ and $\wtilde R_N(t)$. Let $s\equiv s(x,y)\in L^2(\R^d\times\R^d)$ be an integral kernel of operator $S \in \cL^2(\cH)$. Then 
the Wigner function of $S$ is  defined by the formula
\begin{equation}
W_\hb[S](x,\cdot):=\tfrac1{(2\pi)^d}\cF\big(y\mapsto s(x+\tfrac12\hb y,x-\tfrac12\hb y)\big)\quad\text{ for a.e. }x\in\R^d\,,
\end{equation}
where $\cF$ is the Fourier transform on $L^2(\R^d)$.

For each integer $M\ge 1$, we also introduce the dual norm
\begin{equation}\label{eq:dualnorm}
|||f|||_{-M}\!:=\!\sup\left\{\left|\iint_{\R^d\times\R^d}f(x,\xi)\overline{a(x,\xi)}dxd\xi\right|\quad\left|\ba{}&\,\,\,\,\,a\in C_c(\R^d\times\R^d)\,,\,\,\text{ and }
\\
&\max_{|\a|,|\b|\le M\atop |\a|+|\b|>0}\|\d_x^\a\d_\xi^\b a\|_{L^\infty(\R^d\times\R^d)}\!\le\! 1\ea\right.\right\}\,.
\end{equation}

\bigskip
The main results in \cite{golse2019random} is the following theorem.

\begin{theorem}\lb{T-Main}
Assume that $N\ge 2$ and that $V\in C(\R^d)$ is a real-valued function such that
$$
V(z)=V(-z)\text{ for all }z\in\R^d\,,\quad\lim_{|z|\to+\infty}V(z)=0\,,\quad\text{ and }\int_{\R^d}(1+|\om|^2)|\cF (V)(\om)|d\om<\infty\,.
$$
 Let $R_{N,\indc}(t)$ and $\tilde{R}_{N,\indc}(t)$ be the single-particle reduced density operators 
defined in terms of $R_N(t)$ and $R^R_N(t)$ respectively  by \eqref{R1t}.
Then there exists a constant $\g_d>0$ depending only on the dimension $d$ of the configuration space such that, for each $t>0$, one has
\begin{multline}\lb{ErrWig}
|||W_\hb[\mathbb{E} \tilde{R}_{N,\indc}(t)]-W_\hb[R_{N,\indc}(t)]|||_{-[d/2]-3}
\\
\le 2\g_d\Delta te^{6t\max(1,\sqrt{d}L(V))}\L(V)(2+3t\L(V)\max(1,\Delta t)+4\sqrt{d}L(V)t\Delta t)\,,
\end{multline}
where
$$
L(V):=\tfrac1{(2\pi)^d}\int_{\R^d}|\om|^2|\hat V(\om)|d\om\,,\qquad\L(V):=\tfrac1{(2\pi)^d}\int_{\R^d}\sum_{\mu=1}^d|\om^\mu||\hat V(\om)|d\om\,,
$$
with $\omega^\nu$  the $\nu$-th component of $\omega$.
\end{theorem}

\smallskip

This error estimate gives an error {\it independent} of $\hbar$ and $N$. It was also pointed out in \cite{golse2019random} that  the error bound obtained in above theorem is small as $\Delta t\to 0$, \textit{even for moderate values of $N$} for which the factor $\frac1{N-1}$ is insignificant. Therefore the result applies to $N$-body quantum Hamiltonians 
\textit{without the $\frac1{N-1}$ normalization of the interaction potential}, as a simple corollary for each finite value of $N\ge 2$.

\begin{remark}
Note that the dual norm \eqref{eq:dualnorm} is a kind of weak norm. The error bound $\cO(\dt)$ is consistent with the weak error estimate in \cite{jin2021convergence}.
\end{remark}

\subsection{Quantum Monte-Carlo methods}

Computing the ground state energy of a many-body quantum system is a fundamental problem in chemistry. An important tool to determine the ground state energy and electron correlations is the quantum Monte Carlo (QMC) method \cite{von1992quantum,anderson2007quantum}. 

Consider the Hamiltonian,
\begin{equation}\label{eq: ham}
    \cH= \sum_{i=1}^{N} -\frac{\hbar^2}{2m}{\triangle_{x_i}} + \sum_{i\ne j} W(x_i - x_j) + \sum_{i=1}^{N} V_{\mathrm{ext}}(x_i).
\end{equation}
Here  $V_{\mathrm{ext}}$ is the external potential given by
\begin{equation}\label{eq: Vext}
    V_{\mathrm{ext}}(x_i) = \sum_{\alpha=1}^{M} U(x_i - R_\alpha),
\end{equation}
where $R_\alpha$, for instance, can be the position of an atom. 

Up to some global phase factor, the ground state takes real values and is nonnegative everywhere. The ground state and the corresponding eigenvalue can be obtained via the Rayleigh quotient, 
\begin{equation}\label{eq: E-var}
    E= \min_{\Phi_N} \frac{\dsp\int_{(\R^{3})^N} \Phi_N \cH \Phi_N d\underline{x}} 
    {\dsp\int_{(\R^3)^N} |\Phi_N|^2 d\underline{x}},
\end{equation}
where the minimizer $\Phi_N$ corresponds to the ground state  wave function.
The main computational challenge here is the curse of dimensionality due to the high dimensional integral.

In the variational Monte Carlo (VMC) framework, the ground state is approximated by selecting an appropriate ansatz  $\Phi_N \approx \Phi_0$.  Traditionally, $\Phi_0$ is constructed using the one-body wave functions, by taking into  the effect of particle correlations described by the Jastrow factors \cite{foulkes2001quantum}. For example, in the Boson systems like the liquid Helium interacting with a graphite surface \cite{mcmillan1965ground,whitlock1998monte,pang2014diffusion}, the following ansatz has been proven successful,
\begin{equation}\label{eq: ansatz}
    \Phi_0=e^{-J(\underline{x})} \Pi_{i=1}^{N}  \phi(x_i), \quad 
    J(\underline{x})=\frac{1}{2}\sum_{i,j:i\neq j}u(|x_i-x_j|),~~~u(r)= \left( \frac{a}{r} \right)^5 + \frac{b^2}{r^2+c^2}.
\end{equation}
The non-negative one-particle wave function is often taken as
\begin{equation}\label{eq: phii}
\phi(x_i) =  \sum_{\alpha=1}^M e^{-\theta(x_i - R_\alpha)},
\end{equation} 
for some function $\theta$. This form has been used in \cite{whitlock1998monte} and the parameters were obtained by solving a one-dimensional Schr\"odinger equation. With the approximation of $\Phi_N$ being fixed, the multi-dimensional integral is then interpreted as a statistical average. In fact, introducing the probability density function (PDF),
\begin{equation}\label{eq: pdf0}
    p(\underline{x}) \propto |\Phi_0(\underline{x})|^2, 
\end{equation}
the ground state energy is the average of $E_{\text{tot}}$ under $p(\underline{x})$, where
\begin{equation}\label{eq: E0}
     E_\text{tot}(\underline{x})= \frac{\cH \Phi_0} {\Phi_0}.
\end{equation}
Hence,  $E$ can be computed by a Monte Carlo procedure, and such a method is called the VMC, which is a typical QMC method.

In the VMC methods, the ground state is not updated. Instead, one may use another QMC method--the diffusion Monte Carlo (DMC) method  \cite{anderson1975random,reynolds1982fixed}--to compute the ground state and the energy. In particular, one solves a pseudo-time Schr\"odinger equation (TDSE) which is a parabolic equation \cite{reynolds1982fixed}
\begin{equation}\label{eq: tdse}
 \partial_t \Psi_N = (E_T-\cH_N ) \Psi_N. 
\end{equation} 
Here, $t$ represents a fictitious time. The energy shift  $E_T$ is adjusted on-the-fly based on the change of  magnitude of the wave function.  Instead of solving \eqref{eq: tdse} directly, it is often more practical to find $f(\bm r,t)$ with 
\begin{equation}\label{7-33}
    f(\underline{x},t) = \Psi_N(\underline{x},t) \Phi_0(\underline{x}).
\end{equation}
By  choosing  $\Psi_N(\underline{x},0) = \Phi_0(\underline{x})$, $f(\underline{x}, 0)= |\Phi_0|^2 \propto p(\underline{x})$. Hence, a VMC method may be used to initialize $f(\underline{x}, t)$.  Clearly, $f$ solves the following differential equation \cite{reynolds1982fixed},
\begin{equation}\label{eq: f}
 \partial_t f = -\nabla\cdot \big( \frac{\hbar^2}{{m}} v(\underline{x}) f\big) + \frac{\hbar^2}{2m} \nabla^2 f - \big(E_T - {E}_\text{tot}(\underline{x}) \big)  f,
\end{equation}
where $\underline{v}=(v_1, \cdots, v_N)\in \R^{Nd}$ and 
\[
v_i(\underline{x})= \nabla\log\phi(x_i)-\sum_{j: j\neq i} \nabla_{x_i} u(|x_i-x_j|).
\]
The average energy $E(t)$ is then defined as a weighted average,
\begin{equation}
    E(t)= \frac{\dsp\int_{(\R^3)^N} f(\underline{x},t) E_\text{tot}(\underline{x}) d\underline{x}}{\dsp \int_{(\R^3)^N} f(\underline{x}, t) d\underline{x}}, 
\end{equation}
where the correctness can be seen by $E(t)=\int \Psi_N \cH \Phi_0 d\underline{x}/\dsp \int \Psi_N \Phi_0 d\underline{x} $. If $\Psi_N$ is close to the eigenstate, this will be close to $E$.

The key observation is that the dynamics \eqref{eq: f} can be associated with a stochastic process, in which the particles are experiencing birth/death while driven by drift velocity and diffusion.  This process can be implemented by a number of walkers together with birth/death processes  \cite{anderson1975random,reynolds1982fixed}.

\subsubsection{The Random Batch Method for VMC}

With \eqref{eq: ansatz}, the density \eqref{eq: pdf0} can be found as
\begin{equation}\label{eq: pr}
     p(\underline{x}) \propto e^{-2V}, \quad V= - \ln \Psi_0 = -  \sum_i \log \phi (x_i) +  \frac12 \sum_{i} \sum_{j\ne i} u(|x_i - x_j|),
\end{equation}
and the total energy can be  expressed as
\begin{equation}\label{eq: E-tot}
E_\text{tot}(\underline{x})  =  - \frac{\hbar^2}{2m}    \triangle V -  \frac{\hbar^2}{2m}  \| \nabla V\|^2  +   \sum_{i\ne j} W(x_i - x_j) + \sum_{i=1}^{N}   \sum_{\alpha=1}^{M} U(x_i - R_\alpha).
\end{equation}

To sample from $p(\underline{x})$,  one may make use of  the Markov chain Monte Carlo (MCMC) methods. Consider the over-damped Langevin dynamics,
\begin{equation}\label{eq: lgv}
    d {\bm r}_i =  \nabla \log \phi(\bm r_i) dt -  \sum_{j: j\neq i} \nabla_{\bm r_i} u(|\bm r_i - \bm r_j|) dt +  d\bm{W}_i(t), \quad 1 \le i \le N.
\end{equation}
Under suitable conditions \cite{mattingly2002ergodicity}, the dynamical system with potential given by \eqref{eq: pr} is ergodic and the PDF $p(\underline{x})$ in \eqref{eq: pr} is the unique equilibrium measure of \eqref{eq: lgv}. By the classical  Euler-Maruyama method (\cite{kloeden2013numerical}), the underdamped Langevin can be discretized to a Markov Chain:
\begin{equation}\label{eq: EM}
   {\bm r}_i(t+\dt) = {\bm r}_i(t) +    \nabla \log \phi(\bm r_i) \dt -  \sum_{j\ne i} \nabla_{\bm r_i} u(|\bm r_i(t) - \bm r_j(t)|) \dt + \dw_i,  \quad 1 \le i \le N,
\end{equation}
where $\dw_i$ is again sampled from $\mathcal{N}(0, \dt)$. It is clear that  $\mathcal{O}(M+N)$ operations should be taken for each particle at each time step. 

The cost of the above MCMC is high. The strategy in \cite{J-L-QMC} is to apply a RBM strategy with replacement. In particular, at each step, one randomly picks  two particles, $i$ and $j$, and compute their interactions,
 $\nabla_{\bm r_i} u(|\bm r_i - \bm r_j|)$,
 then updates their positions as follows, 
\begin{equation}
\label{eq-1}
\left\{
\begin{aligned}
    {\bm r}_i(t+\dt) &= {\bm r}_i(t) + \nabla \log \phi(\bm r_i) \dt  + (N-1) \nabla_{\bm r_i} u(|\bm r_i - \bm r_j|)  \dt +  \dw_i, \\
    {\bm r}_j(t+\dt) &= {\bm r}_j(t) + \nabla \log \phi(\bm r_j) \dt  + (N-1) \nabla_{\bm r_j} u(|\bm r_i - \bm r_j|)  \dt +  \dw_j.
\end{aligned}
\right.
\end{equation}
For the one-body term $\nabla\log \phi(\bm r_i)$,
\begin{equation}\label{eq: 1-body}
   \nabla \log \phi(\bm r_i)  = \dsp \sum_{\alpha=1}^M - \nabla \theta(\bm r_i - R_\alpha) q_\alpha^i, \quad q_\alpha^i = \frac{e^{-\theta(\bm r_i - R_\alpha)}}{\sum_{\beta=1}^M e^{-\theta(\bm r_i - R_\beta)}},
\end{equation}
where the coefficients $q_\alpha^i$'s are non-negative and $\sum_\alpha q_\alpha^i =1$. To reduce the cost, one may further use a direct Monte-Carlo method: pick {\it just one} term $\alpha$ randomly. Specifically, assume that one starts with $\alpha$ and  computes
$e_{\mathrm{old}}=\theta(\bm r_i - R_\alpha)$, and then one randomly picks $1 \le \beta \le M$, and computes
$e_{\mathrm{new}}=\theta(\bm r_i - R_\beta)$.  $\beta$ is accepted with probability 
\begin{equation}\label{eq:pacc}
p_{\mathrm{acc}}\propto \exp\big[ -( e_{\mathrm{new}}-e_{\mathrm{old}}) \big]. 
\end{equation}
\medskip
For the detailed algorithm see \cite{J-L-QMC}. As a result of the random sampling of the one- and two-body interactions, updating the position of each particle {\it only requires $\mathcal{O}(1)$ operations} per time-step.  Another practical issue emerges when the interaction $u(|x|) $ has a singularity near zero. One can use the splitting idea as mentioned in section \ref{subsec:splitting}, i.e., applying RBM only to the long-range smooth part.

It was shown  in \cite{J-L-QMC} that the above random batch algorithm, when applied to one batch of two particles, has the same accuracy as the Euler-Maruyama method over a time step of 
$2\dt/N.$  One full time step in Euler-Maruyama method corresponds to $N/2$ such steps in the random batch algorithm. This corresponds to the random batch method with replacement. 

We show a  numerical experiment performed in \cite{J-L-QMC} on ${}^4$He atoms interacting with a two-dimensional lattice.  The CPU times taken to move the 300 Markov chains for 1000 steps were compared.  In this comparison, the cost associated with the energy calculations was excluded in the random batch and Euler-Maruyama methods. From Table \ref{vmc}, one clearly sees that the RBM is more efficient than the Euler-Maruyama method. It is much more efficient than the random walk Metropolis-Hastings algorithm, mainly because the latter method requires the calculation of the energy at {\it every } step. 
\begin{table}[thp]

\caption{CPU times (seconds) for several VMC methods. \label{vmc} }
\begin{tabularx}{\textwidth}{b|s|s|s}
\hline\hline
  & Random Walk Metropolis-Hastings  & Euler-Maruyama &   Random Batch   \\
 \hline
 CPU time for a 1000-step sampling period  &  1503  & 469 &  54  \\
  \hline\hline
 \end{tabularx}
 \end{table}

\subsubsection{The Random Batch Method for DMC}

Viewing \eqref{eq: f}, one may consider an ensemble of $L$ copies of the system, also known as walkers \cite{anderson1975random}. For each realization, one first solves the SDEs corresponding to the drift and diffusion, which
is the same as the overdamped Langevin as in VMC  up to a time scaling. Hence, the same Random Batch Algorithm in the VMC can be used for this part.

The relaxation term $-(E_T-E_{\mathrm{tot}})f$ is then done by using a birth/death process to determine whether a realization should be removed or duplicated.  For each walker, one computes a weight factor,
\begin{equation}
 w(t+\Delta t) = \exp \left[ \Delta t\big(E_T -  \tfrac{1}{2}({E}_\text{tot}(\bm r) +{E}_\text{tot}(\bm r'))  \big) \right].
\end{equation}
This weight determines how the walker should be removed or duplicated. See \cite{J-L-QMC} for more details.
The primary challenge is that computing the energy at each step requires $\mathcal{O}((N+M)N)$ operations in order to update the position of $N$ particles. To reduce this part of the computation cost,  one rewrites the total energy as
 \begin{equation}\label{eq: e-parti}
 {E}_\text{tot}(\bm r)= \sum_{i=1}^N E_1(\bm r_i) + \sum_{1\le i<j\le N} E_2(\bm r_i, \bm r_j) + \sum_{1\le i <j<k\le N} E_3(\bm r_i, \bm r_j, \bm r_k),
\end{equation}
where
\begin{equation}
\begin{split}
& E_1(\bm r_i)= - \frac{\hbar^2}{2m} \nabla^2 \ln \phi(\bm r_i) - \frac{\hbar^2}{2m} |\nabla  \ln \phi(\bm r_i) |^2 +  \sum_{\alpha=1}^{M} U(\bm r_i - R_\alpha),\\
& E_2(\bm r_i, \bm r_j)= - \frac{\hbar^2}{m} \nabla^2 \ln u( r_{ij}) + \frac{\hbar^2}{m} \big(\nabla  \ln \phi(\bm r_i)-\nabla  \ln \phi(\bm r_j)\big) \cdot \nabla  u( r_{ij})  + \frac{\hbar^2}{m}  |\nabla  u( r_{ij})|^2 + W( r_{ij}),\\
& E_3(\bm r_i, \bm r_j, \bm r_k)=  \frac{\hbar^2}{m} \Big[ \nabla  u( r_{ij}) \cdot \nabla  u( r_{ik}) + \nabla  u(r_{ji}) \cdot \nabla  u( r_{jk}) + \nabla  u( r_{ki}) \cdot \nabla  u( r_{kj}) \Big].
\end{split}
\end{equation}
Here,  $\bm r_{ij}=\bm r_i -\bm r_j$ and $r_{ij}=|\bm r_{ij}|$.
The three-body terms arise because of the $\|\nabla V\|^2$ term in \eqref{eq: E-tot}.

 In the random batch algorithm proposed in \cite{J-L-QMC}, one randomly picks a batch $C_I$ with three particles:  $C_I=\{ i, j, k\}.$  One first updates the position of the three particles (drift and diffusion) by solving the overdamped Langevin dynamics using the random batch algorithm with batch size 3.  Then,  one then defines a {\it local} energy,
\begin{equation}\label{eq: EIijk}
\begin{aligned}
  {E}_I(\bm r_i, \bm r_j, \bm r_k) =&  E_1(\bm r_i) + E_1(\bm r_j) + E_1(\bm r_k)  \\
          & + \tfrac{N-1}{2} \Big[ E_2(\bm r_i,\bm r_j) + E_2(\bm r_j,\bm r_k) + E_2(\bm r_k,\bm r_i) \Big],\\
          & + \tfrac{(N-1)(N-2)}{2} E_3(\bm r_i,\bm r_j,\bm r_k).   
\end{aligned}
\end{equation} 
where in $E_1$, the sum $\sum_{\alpha=1}^M$ can be further reduced by a 
mini-batch strategy.  Computing this local energy is clearly $\cO(1)$.
To avoid frequent removal and duplication of walkers,  the branching process is applied  after $N/3$ batches of particles are updated. In this case, the weight function is defined by collecting the local energy from each batch (denoted by $I_m$ here),
\begin{equation}
    w(\bm r) = \exp \left[\Delta t \big( E_T  - \widetilde{E}_\text{tot} \big)\right], \quad \widetilde E_\text{tot}=\sum_{m=1}^{N/3} {E}_{I_m}.
\end{equation}
Because of the smallness of $\dt$, the expectation of $w_I$ equals $w(t+\dt)$ modulus an error of  $\cO(\dt^2)$. 
See \cite{J-L-QMC} for the verification using the Green's functions.  

\begin{figure}[htbp]
\centering
    \includegraphics[scale=0.12]{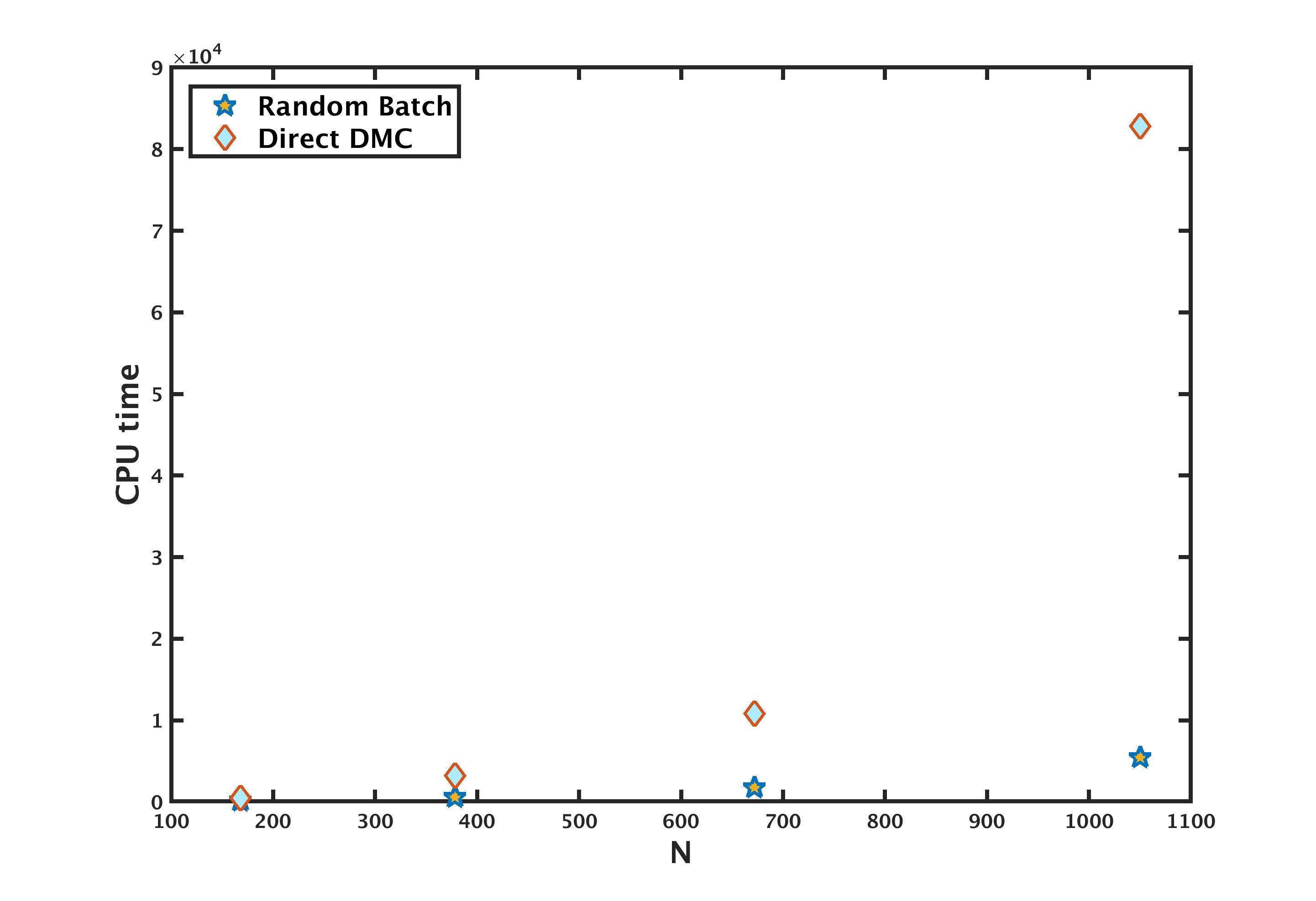}
    \caption{A comparison of the CPU runtime (in seconds) for running 1000 steps of DMC. }
    \label{fig:CPU}
\end{figure}

The detailed algorithm can be found in \cite{J-L-QMC} and we omit it here. 
Now we show a  test of the RBM-DMC algorithm conducted in \cite{J-L-QMC}, which compares the results with the direct DMC method. For the initialization,  a VMC method using the ansatz \eqref{eq: ansatz} for the wave function $\Phi_0$ was first applied.  The random walk Metropolis-Hastings Monte Carlo method is used in both methods so that they start at the same states. 300 ensembles are created by sub-sampling one sample out of every 500 steps from the VMC runs to avoid correlations among the ensembles. For both methods,  $\dt=10^{-4}$ was used and  $200,000$ steps of  simulations were run.   
The CPU run-time is recorded for various system sizes. More specifically,  the system size is increased from the original 168 particles, to $N= 378$,  $N=672$ and $N=1050$ particles, and in each case,  the direct DMC and the RBM-DMC were run for 1000 steps. As shown in Figure \ref{fig:CPU}, the CPU time for the direct DMC method increases much more rapidly as $N$ increases. 

\section*{Acknowledgement}
S. Jin was partially supported by the NSFC grant No.12031013.  The work of L. Li was partially sponsored by NSFC 11901389, 11971314, and Shanghai Sailing Program 19YF1421300. Both authors were also supported by  Shanghai Science and Technology Commission Grant No. 20JC144100.

\bibliographystyle{plain}
\bibliography{sdealg}

\end{document}